\documentclass[12pt,letterpaper,reqno]{amsart}
\usepackage{amsthm,amssymb,amsmath,amscd}
\usepackage{graphicx} 
\usepackage[x11names]{xcolor}
\usepackage{tikz}
\usepackage{fourier}
\usepackage{subfigure}
\usepackage{tikz}
\usepackage{overpic}
\usepackage[breaklinks,colorlinks,citecolor=blue,linkcolor=blue,
 urlcolor=teal]{hyperref} 
 \usepackage{thmtools}
\colorlet{mycolor}{DarkOrchid3}
\usetikzlibrary{ decorations.markings}
\usetikzlibrary{arrows,arrows.meta,calc,decorations.markings} 
\theoremstyle{plain}
\newtheorem{theorem}{Theorem}[section]
\newtheorem{proposition}[theorem]{Proposition}
\newtheorem{lemma}[theorem]{Lemma}
\newtheorem{corollary}[theorem]{Corollary}
\newtheorem{question}[theorem]{Question}

\theoremstyle{definition}
\newtheorem{definition}[theorem]{Definition}
\newtheorem{remark}[theorem]{Remark}

\newtheorem*{ack}{Acknowledgments}
\newtheorem*{ootp}{Outline of the Paper}

\newcommand{\co}{\colon \thinspace}

\newcommand{\abs}[1]{\lvert {#1} \rvert}

\renewcommand{\leq}{\leqslant}
\renewcommand{\geq}{\geqslant}
\renewcommand{\epsilon}{\varepsilon}
\renewcommand{\phi}{\varphi}

\newcommand{\N}{{\mathbb N}}
\newcommand{\Q}{{\mathbb Q}}
\newcommand{\R}{{\mathbb R}}
\newcommand{\C}{{\mathbb C}}
\newcommand{\Z}{{\mathbb Z}}
\newcommand{\M}{{\mathcal M}}
\newcommand{\Hy}{{\mathbb H}}
\DeclareMathOperator{\Isom}{Isom}
\DeclareMathOperator{\Aut}{Aut}

\numberwithin{theorem}{section}
\numberwithin{equation}{section}

\begin{document}

\title[Minimal complexity $3$-manifolds]{Minimal complexity cusped hyperbolic $3$-manifolds \\ with geodesic boundary}

\author{Anuradha Ekanayake}
\author{Max Forester}
\author{Nicholas Miller}
\address{Department of Mathematics\\University of Oklahoma\\Norman, OK 73019}
\email{anuradha.ekanayake@ou.edu, mf@ou.edu, nickmbmiller@ou.edu}


\begin{abstract}
In the early 2000s, Frigerio, Martelli, and Petronio studied $3$-manifolds of smallest combinatorial complexity that admit hyperbolic structures. As part of this work they defined and studied the class $\M_{g,k}$ of smallest complexity manifolds having $k$ torus cusps and connected totally geodesic boundary a surface of genus $g$. In this paper, we provide a complete classification of the manifolds in $\M_{k,k}$ and $\M_{k+1,k}$, which are the cases when the genus $g$ is as small as possible.
In addition to classifying manifolds in $\M_{k,k}$, $\M_{k+1,k}$, we describe their isometry groups as well as a relationship between these two sets via Dehn filling on small slopes. 
Finally, we give a description of important commensurability invariants of the manifolds in $\M_{k,k}$.
\end{abstract}


\maketitle

\thispagestyle{empty}


\section{Introduction}


In the study of hyperbolic $3$-manifolds, it is of interest to find those which are minimal with respect to some natural measure of complexity. 
There is a lineage of results where this complexity is taken to be volume; for example Gabai, Meyerhoff, and Milley \cite{GMM} showed that the Weeks manifold is the smallest volume closed hyperbolic $3$-manifold and Cao and Meyerhoff \cite{CaoMeyerhoff} showed that the figure-eight knot complement and its sister are the minimal volume $1$-cusped hyperbolic $3$-manifolds.
In the sequence of papers \cite{fmp-small,fmp-complexity,fp} as well as other related works, Frigerio, Martelli, and Petronio  attempted to classify manifolds which have minimal complexity in a certain combinatorial sense described below; informally, they are minimally triangulated with respect to number of tetrahedra.
For hyperbolic $3$-manifolds with totally geodesic boundary, these notions of complexity appear to be intimately related, as has been shown in small volume examples by Kojima and Miyamoto \cite{KojimaMiyamoto} and the first author \cite{Anuthesis}.
Indeed, in these small volume examples, the manifolds which minimize volume as a notion of complexity also minimize combinatorial complexity and vice versa.
This paper continues the study of minimal combinatorial complexity manifolds, which we now describe more precisely.

Throughout, we use $\Delta$ to denote the standard tetrahedron and $\mathring\Delta$ to denote $\Delta$ with its vertices removed.
Then an \emph{ideal triangulation}, $\mathcal{T}$, of a compact $3$-manifold $M$ with non-empty boundary is a gluing of finitely many $\mathring\Delta$ via simplicial face pairings whose union is the interior of $M$.
Let $|\mathcal{T}|$ denote the number of $\mathring\Delta$ in a triangulation $\mathcal{T}$, then we define the combinatorial complexity of a manifold $M$ by the formula
$$c(M)=\min\left\{~\!|\mathcal{T}|~\middle|~\mathcal{T}\text{ ideally triangulates }M\right\}.$$
Throughout this paper, we let $\M^{cco}$ denote the set of homeomorphism types of compact, connected, orientable $3$-manifolds (possibly with boundary).

In \cite[\S2]{fmp-dehn}, it was shown that $g+k$ is the minimal size of an ideal triangulation of a manifold $M\in\M^{cco}$ with the property that $\partial M=\Sigma_g\sqcup\left(\sqcup_{i=1}^k\mathbb{T}^2\right)$, where $\Sigma_g$ is a surface of genus $g\ge 2$, $\mathbb{T}^2$ is the $2$-torus, and $k$ is any non-negative integer.
Consequently, the lower bound $c(M)\ge g+k$ always holds for such manifolds.
Taking this into account, one defines the set of minimal complexity $3$-manifolds with boundary a union of $k$ tori and a genus $g$ surface by
$$\mathcal{M}_{g,k}=\left\{[M]\in\M^{cco}~\middle|~ \partial M=\Sigma_g\sqcup\left(\sqcup_{i=1}^k\mathbb{T}^2\right),~c(M)=g+k \right\}.$$
In the sequel, we frequently blur the distinction between $M$ and its homeomorphism class when it is clear from context.

The connection between manifolds $M\in\M_{g,k}$ and hyperbolic $3$-manifolds is given by \cite[Thm 1.2(1)]{fmp-dehn}, relying on the work of \cite{fp}, wherein Frigerio, Martelli, and Petronio show that such a combinatorial ideal triangulation gives rise to a geometric triangulation by geometric partially truncated tetrahedra, thereby endowing $M$ with the structure of a hyperbolic $3$-manifold with $k$ cusps and totally geodesic boundary of genus $g$.
This hyperbolic structure is then unique by Mostow rigidity (see Section \ref{subsection:tetrabackground}).
Moreover, they show in \cite[Thm 1.2(2)]{fmp-dehn} that their procedure is canonical in the sense that it yields the canonical Kojima decomposition \cite{kojima}, a generalization of the Epstein--Penner decomposition of a finite-volume hyperbolic $3$-manifold.

Fixing a natural number $C\ge 2$, in \cite{fmp-dehn} the aforementioned authors studied the growth rate of the number of hyperbolic $3$-manifolds of complexity $C$, culminating in \cite[Cor 1.5]{fmp-dehn}, where they show that the growth rate of such manifolds is $C^C$.
The proof of this follows from \cite[Thm 1.4]{fmp-dehn}, where it is shown that for any fixed non-negative integer $k$, the growth type of $|\M_{g,k}|$ is $g^g$.
As $\M_{g,k}$ is a subset of the manifolds of combinatorial complexity $g+k$, this immediately yields the growth rate with respect to complexity.

As opposed to understanding growth rates of all hyperbolic $3$-manifolds with complexity $C$, in this paper we are interested in the extremal manifolds $M$ of a fixed complexity $C$. 
Here we use the term extremal to mean that either $M$ has the least number of cusps possible for a manifold of complexity $C$ or the largest number of cusps possible for a manifold of complexity $C$.
The case where $M$ has the least number of cusps possible is already resolved by the work of \cite[Thm 1.4]{fmp-dehn} upon setting $k=0$.
As previously mentioned, these manifolds have a growth rate of $C^C$ and there are a plethora of such manifolds as $g$ increases.

In contrast, in this paper we show that there are relatively few homeomorphism types of hyperbolic $3$-manifolds having the largest number of cusps possible for a manifold of complexity $C$ (equivalently, with the smallest possible genus of the totally geodesic boundary).
Such manifolds always lie in $\M_{k,k}$ if $C=2k$ is even or in $\M_{k+1,k}$ if $C=2k+1$ is odd.
With this in mind, we now describe our main results.

In the case where the complexity is even, the aforementioned manifolds will lie in $\M_{k,k}$.
For these manifolds, it was already shown in \cite[Prop 1.3]{fmp-dehn} that $\M_{k,k}=\emptyset$ if and only if $k$ is odd.
Our first result gives a complete description of the manifolds in $\M_{k,k}$ for $k$ even, which is the content of the following.

\begin{restatable}{theorem}{thmkkeven}\label{thm-kkeven}
If $k$ is even then $\abs{\M_{k,k}} = 1$. That is, there is a unique homeomorphism type of manifold in $\M_{k,k}$.
\end{restatable}

The method of proof will allow us to also provide a complete description of their isometry groups.
In what follows, let $D_n$ denote the dihedral group of order $2n$. 

\begin{restatable}{theorem}{thmisometrieskk}\label{thm-isometries-kk}
Let $M$ be the unique manifold in $\M_{k,k}$ for $k$ even. The isometry group of $M$ is isomorphic to $D_{3k}$, with orientation-preserving subgroup $D_{3k/2}$.  
\end{restatable}
 
Next we address the case where the complexity is odd so that the manifolds lie in $\M_{k+1,k}$.
In contrast to the situation for $\M_{k,k}$, we show that manifolds in $\M_{k+1,k}$ exist for any $k\in\N$ and moreover give a complete classification of all possible topological types.

\begin{theorem}\label{thm-kplusone}
For any $M \in \M_{k+1,k}$ there is a canonical partition of the set of torus components of $\partial M$ into two sets (one of which may be empty). Let $i(M)$ and $j(M)$ denote the cardinalities of these two sets, where $i(M)\leq j(M)$. Then
\begin{enumerate}
\item $i(M)$ and $j(M)$ are not both even; 
\item manifolds $M_1, M_2 \in \M_{k+1,k}$ are homeomorphic if and only if $(i(M_1),j(M_1)) = (i(M_2),j(M_2))$;  
\item every pair $(i,j)$ with $k = i+j$, $0 \leq i \leq j$, and $i,j$ not both even is $(i(M),j(M))$ for some $M \in \M_{k+1,k}$. 
\end{enumerate}
\end{theorem}

In particular, counting the number of all such pairs in Theorem \ref{thm-kplusone} immediately yields the following corollary. 

\begin{corollary}\label{cor:kplusonecount}
  For any $k \geq 1$, we have that
  \[ \abs{\M_{k+1,k}} \ = \ \begin{cases}
    (k+1)/2, & \text{ if } \ k \text{ is odd}, \\
    k/4, & \text{ if } \ k \equiv 0 \pmod 4, \\
    (k+2)/4, & \text{ if } \ k \equiv 2 \pmod 4. \\
    \end{cases}\]
\end{corollary}

For context, Corollary \ref{cor:kplusonecount} in the case that $k=1$ was previously shown computationally in work of Frigerio, Martelli, and Petronio \cite{fmp-small}.
To our knowledge the computation of these sets is new for all $k$ strictly bigger than $1$.

Similar to the above, we are also able to determine the isometry groups of manifolds in $\M_{k+1,k}$ as a consequence of the methods developed in Section \ref{sec:kplusone} (see Corollary \ref{cor:isometrieskplusone}).

\begin{theorem}\label{thm:koplusoneisom}
For every $M\in\mathcal{M}_{k+1,k}$, we have that
$$\Isom(M)\cong\begin{cases}
\{ 1 \}, & \text{ if } \ i(M)\neq j(M),\\
\Z/2\Z, & \text{ if } \ i(M)=j(M),
\end{cases},$$
where this isometry group is with respect to the unique hyperbolic structure on $M$.
\end{theorem}

Next, we consider the relationship between manifolds in $\M_{k,k}$ and those in $\M_{k,k-1}$. Provided $k$ is even, it was shown in \cite[Thm 1.7]{fmp-dehn} that if one performs Dehn filling on one cusp of $M \in \M_{k,k}$ with slope $p/q \in \{-2,-1/2,1/3,2/3,3/2,3\}$, then the resulting manifold is always contained in $\M_{k,k-1}$. In this paper, we prove that for each of these six slopes, the resulting filled manifold is always the same and identify its homeomorphism type amongst manifolds in $\M_{k,k-1}$. 


\begin{restatable}{theorem}{dehnfill}\label{thm-dehnfill}
Let $k$ be an even natural number and $M\in\mathcal{M}_{k,k}$ the unique manifold furnished by Theorem \ref{thm-kkeven}. Then for all $p/q\in \{-2,-1/2,1/3,2/3,3/2,3\}$, all manifolds obtained by $p/q$-Dehn filling on any single cusp of $M$ are homeomorphic to the same manifold $M_{0}$ in $\M_{k,k-1}$.
Moreover, this manifold always corresponds to the trivial partition of its set of torus boundary components, that is, the unique manifold with invariants $i(M_0)=0$, $j(M_0)=k-1$ in the description of Theorem \ref{thm-kplusone}.
\end{restatable}


Note that the special case that $k=2$ is immediate, since $|\mathcal{M}_{2,2}|=|\mathcal{M}_{2,1}|=1$, however in all other cases Theorem \ref{thm-dehnfill} is new.

Finally, we give an understanding of the topological invariants associated to the hyperbolic structure of the unique manifold in $\M_{k,k}$ for $k$ even. This is the content of the following theorem, where we refer the reader to Section \ref{sec:examples} for the relevant definitions.

\begin{restatable}{theorem}{kevenquasiarithmetic}\label{thm:kevenquasiarithmetic}
Suppose that $k$ is even and $M$ is the unique element of $\mathcal{M}_{k,k}$.
Then the invariant trace field of $M$ is $L_k=\Q(\sqrt{-3},\cos\left(\frac{2\pi}{3k}\right))$ and the invariant quaternion algebra of $M$ is $\mathrm{Mat}_2(L)$. 
Moreover, $M$ has non-integral traces if and only if $k$ is a power of\/ $2$.
In particular, $M$ is quasi-arithmetic if and only if $k=2$ and is non-arithmetic for all values of\/ $k$.
\end{restatable}

Combining the work from \cite[\S 2]{fmp-dehn}, \cite{fp}, and a volume formula due to Ushijima \cite{Ushijima}, one can compute that the unique manifold $M\in\M_{k,k}$ has volume $V_k$, whose explicit expression can be written using dilogarithms and is computed explicitly in Equation \eqref{eqn:volumeMkk} in Section \ref{subsection:volume}. 
In the same section, we also show that this volume has linear growth rate with respect to $k$.

Based on the work done in this paper and its seeming connection to volume minimization, it seems entirely possible that the manifolds in $\M_{k,k}$ and $\M_{k+1,k}$ should also serve as volume minimizers.
We ask this formally in the following question.

\begin{question}\label{ques:minimalvol}
Let $M$ be the unique manifold in $M_{k,k}$, is $M$ minimal volume amongst hyperbolic $3$-manifolds with $k$-cusps and connected, totally geodesic boundary? Is it the unique one?
Similarly, are any of the manifolds in $M_{k+1,k}$ minimal volume in the class of hyperbolic $3$-manifold with $k$-cusps and a connected, totally geodesic boundary? 
\end{question}

It is known by work of Kojima and Miyamoto \cite{KojimaMiyamoto}, building on work of Fujii \cite{Fujii}, that there are eight manifolds in $\M_{2,0}$ that are minimal amongst all compact hyperbolic $3$-manifolds with totally geodesic boundary and by work of the first author \cite{Anuthesis} that the unique manifold in $\M_{2,1}$ is minimal amongst all $1$-cusped manifolds hyperbolic $3$-manifolds with totally geodesic boundary.
We are therefore tempted to conjecture an affirmative answer to Question \ref{ques:minimalvol}; however, based on the lack of numerical examples in this volume regime, we prefer to state it as a question.

\begin{ootp}
In Section \ref{sec:background}, we will describe the necessary background for the rest of the paper. Of particular importance will be the combinatorial description of triangulations of manifolds from $\M_{g,k}$ in Section \ref{subsection:tetrabackground} and an explicit description of the spines of these manifolds in Section \ref{subsection:spinebackground}. 
In Section \ref{sec:kkeven}, we will prove Theorem \ref{thm-kkeven} by giving an explicit combinatorial description of manifolds in $\M_{k,k}$ and Theorem \ref{thm-isometries-kk} using an explicit description of the big cell of the spine of these manifolds.
Using similar techniques, in Section \ref{sec:kplusone} we will prove Theorem \ref{thm-kplusone} and Theorem \ref{cor:isometrieskplusone} by analyzing the possible combinatorial triangulations of manifolds in $M_{k+1,k}$ as well as their corresponding big cells.
In Section \ref{sec:dehnfilling}, we prove Theorem \ref{thm-dehnfill} by analyzing how the aforementioned spine changes under Dehn filling on small slopes.
Finally, in Section \ref{sec:examples} we will prove Theorem \ref{thm:kevenquasiarithmetic} by analyzing the corresponding invariants for certain reflection groups which are commensurable to the manifold in $\M_{k,k}$.
\end{ootp}

\begin{ack}
The third author is partially supported by NSF DMS--2405264.
\end{ack}


\section{Background}\label{sec:background}


\subsection{Generalities on ideal triangulations of hyperbolic manifolds}\label{subsection:tetrabackground}

In this section we recount work of Frigerio, Martelli, and Petronio \cite{fmp-dehn}, building on work of Kojima \cite{kojima} and Frigerio and Petronio \cite{fp}, which describes the combinatorics of minimal complexity ideal triangulations of manifolds in $\mathcal{M}_{g,k}$ as well as a procedure for associating a hyperbolic structure to such a triangulation.
In what follows, we retain the notation and definitions from the introduction.

Let $M\in\mathcal{M}_{g,k}$ and $\mathcal{T}$ a corresponding minimal ideal triangulation. Using the fact that $|\mathcal{T}|=g+k$, Frigerio, Martelli, and Petronio obtain an explicit combinatorial description of $\mathcal{T}$.
This culminates in the following proposition (see \cite[Prop 2.2]{fmp-dehn}), where we remind the reader that the incidence number of an edge is the number of tetrahedra incident to it, counted with multiplicity.


\begin{proposition}\label{prop-combinatorics}
Suppose $M\in\mathcal{M}_{g,k}$ and let $\mathcal{T}$ be a corresponding minimal ideal triangulation. Then
\begin{enumerate}
\item For each $i=1,\dots, k$ there are exactly two tetrahedra of $\mathcal{T}$ with $3$ vertices on $\Sigma_g$ and one vertex on the $i^{th}$ copy of $\mathbb{T}^2$; the remaining $g-k$ tetrahedra have all of their vertices on $\Sigma_g$;
\item The triangulation $\mathcal{T}$ has precisely $k+1$ edges $e_0,\dots, e_k$. Moreover, $e_0$ has both endpoints on $\Sigma_g$ with incidence number $6g$ and each $e_i$ for $1\le i\le k$ has one vertex on $\Sigma_g$, one vertex on the $i^{th}$ copy of $\mathbb{T}^2$, and incidence number $6$.
\end{enumerate}
\end{proposition}


Given a manifold $M\in\mathcal{M}_{g,k}$ with corresponding ideal triangulation $\mathcal{T}$, there is a procedure for geometrizing each tetrahedron so as to give $M$ a hyperbolic structure. We recall this procedure below, first describing how to endow an individual tetrahedron $\Delta$ with a hyperbolic structure.


\begin{figure}[t]
\centering
\begin{subfigure}
\centering
\includegraphics[width=1.75in]{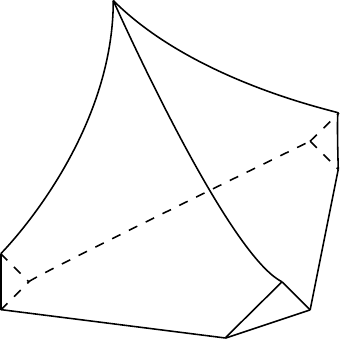}
\end{subfigure}
\hspace{.75in}
\begin{subfigure}
\centering
\includegraphics[width=1.75in]{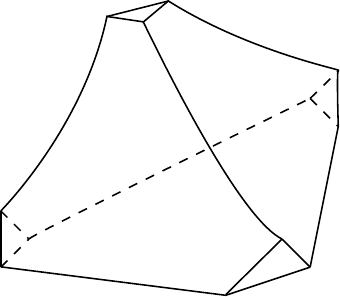}
\end{subfigure}
\caption{Two examples of partially truncated tetrahedra. The one on the left has a unique ideal vertex and the one on the right is a compact tetrahedron}
\label{fig-pttexamples}
\end{figure}


Let $\Delta$ be as above and let $\mathcal{I}$ denote a subset of the vertex set $V(\Delta)$ of $\Delta$, which we will refer to as the \emph{ideal vertices of $\Delta$}. 
Then we obtain a subset $\Delta^*$ of $\Delta$ by 1) deleting each vertex in $\mathcal{I}$ from $\Delta$ and 2) deleting a small open neighborhood of each vertex $v\in V(\Delta)\setminus \mathcal{I}$, creating a new triangular face; see Figure \ref{fig-pttexamples}.
We call $\Delta^*$ a \emph{partially truncated tetrahedron} and omit reference to the set $\mathcal{I}$ when it is clear from context.
When $\mathcal{I}=\emptyset$, we will call $\Delta^*$ a \emph{compact tetrahedron} and otherwise refer to $\Delta^*$ as a \emph{non-compact tetrahedron}.

The faces of $\Delta^*$ fall into two types: the \emph{triangular faces} corresponding to the vertices of $\Delta$ not contained in $\mathcal{I}$, and the \emph{original faces}, which are either hexagons or polygons with some ideal vertices. The edges of $\Delta^*$ are either \emph{boundary edges} (the edges of the triangular faces) or \emph{original edges}.

A \emph{geometric realization of $\Delta^*$} will refer to an embedding of $\Delta^*$ into $\Hy^3$ such that each face of $\Delta^*$ is totally geodesically embedded and each boundary edge has dihedral angle $\pi/2$. Important for us will be the following two facts, which follow from \cite{fp}:
\begin{enumerate}
\item If $|\mathcal{I}|=1$ then, up to isometry of $\Hy^3$, there is a unique geometric realization of $\Delta^*$ for which every original edge whose closure contains the unique ideal vertex has dihedral angle $\pi/3$. This is the only type of non-compact tetrahedron that will arise in our constructions in the sequel.
\item If $\Delta^*$ is a compact tetrahedron then, up to isometry of $\Hy^3$, there is a unique geometric realization where each original edge has dihedral angle $\alpha$ for some $\alpha\in(0,\pi/3)$.
\end{enumerate}
We will also occasionally abusively use the terminology describing faces and edges of $\Delta^*$ for the tetrahedron $\Delta$ as well as, that is, to describe the similar edges/faces before taking the geometric realization.

Using facts (1) and (2) and the combinatorial description of an ideal triangulation with minimal complexity from Proposition \ref{prop-combinatorics}, Frigerio, Martelli, and Petronio prove the following theorem, which is the content of \cite[Thm 1.2(1)]{fmp-dehn}.

\begin{theorem}\label{thm-hypvol}
Suppose that $M\in\mathcal{M}_{g,k}$. Then $M$ has a unique complete hyperbolic structure realizing $M$ as a hyperbolic $3$-manifold with $k$ cusps and a totally geodesic boundary of genus $g$.
Moreover, the volume of $M$ in this hyperbolic structure depends only on $g$ and $k$.
\end{theorem}

In \cite{kojima}, Kojima proved that any complete hyperbolic manifold with non-empty totally geodesic boundary admits a canonical decomposition into geometric partially truncated polyhedra, which we call the \emph{canonical Kojima decomposition}. 
We will not recount the details of this decomposition herein, as we will not need them in the sequel, however important to us will be the following result of Frigerio, Martelli, and Petronio \cite[Thm 1.2(2)]{fmp-dehn}. 


\begin{theorem}\label{thm-kojima}
Suppose that $M\in\mathcal{M}_{g,k}$. Then $M$ has a unique ideal triangulation $\mathcal{T}$ with $g+k$ tetrahedra.
Moreover, the geometric realization of $\mathcal{T}$ gives the canonical Kojima decomposition of $M$.
\end{theorem}


In particular, one obtains the following corollary, which allows us to reduce the study of the homeomorphism types of a manifold $M$ in $\mathcal{M}_{g,k}$ to the question of combinatorial isomorphism.

\begin{corollary}\label{cor:combinatorialiso}
Suppose that $M$, $M'$ have triangulations $\mathcal{T}$, $\mathcal{T}'$ (respectively) as in Theorem \ref{thm-kojima}. 
Then $M$ and $M'$ are homeomorphic if and only if the triangulations $\mathcal{T}$ and $\mathcal{T}'$ are combinatorially equivalent.
\end{corollary}


\subsection{Generalities on spines}\label{subsection:spinebackground}

Let $\Delta^*$ be a partially truncated tetrahedron. In this section we describe how to endow $\Delta^*$ with a dual combinatorial structure that we call a vertebra.
The union of all vertebrae arising from a given ideal triangulation will be a subset called the spine, which is a subpolyhedral complex onto which the manifold collapses. 
We now describe this construction in more detail.

\begin{figure}[t]
\centering
\includegraphics[width=4in, height=2in]{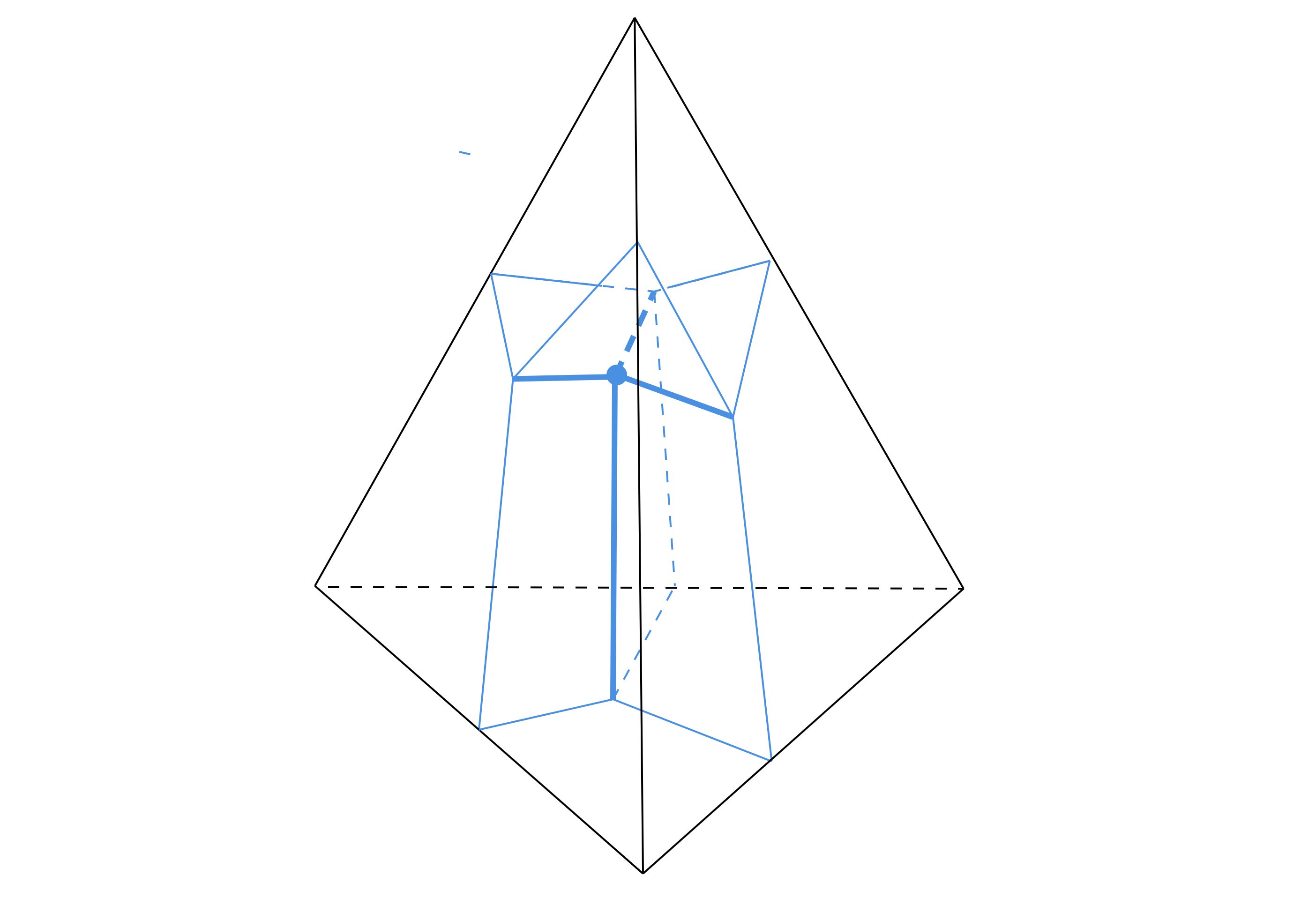}
\caption{A tetrahedron outlined in black with its associated spinal vertex, four spinal half edges, and six diamonds all drawn in blue.}
\label{fig-vertebra}
\end{figure}

Momentarily considering $\Delta$ with its Euclidean structure, we call the barycenter of $\Delta^*$ the \emph{spinal vertex}. 
One can then construct the \emph{spinal half-edges}, each of which is the convex hull of the spinal vertex and the barycenter of a face. There are four spinal half-edges.
Next there are $2$-dimensional subsets called \emph{diamonds}: for each edge $e$ of $\Delta$, the diamond \emph{dual to $e$} or \emph{meeting $e$} is the convex hull of the set consisting of the midpoint of $e$ and the two spinal half-edges meeting the faces adjacent to $e$. There are six diamonds. 
The \emph{vertebra} of $\Delta$ is the union of the spinal vertex, the four spinal half-edges, and the six diamonds. 
See Figure \ref{fig-vertebra} for a visual representation.


Following Matveev's theory of spines \cite{matveev1990complexity}, a \emph{spine} of a manifold is a subpolyhedron of the manifold onto which the manifold collapses. 
For manifolds $M\in\mathcal{M}_{g,k}$ with ideal triangulation $\mathcal{T}$, the union of the vertebrae of the $g+k$ partially truncated tetrahedra in $\mathcal{T}$ is a spine $P$ of $M$ \cite[\S 3]{fmp-dehn}. It has a cell structure which is dual to $\mathcal{T}$ in the following sense. For each edge $e$ of $\mathcal{T}$, the diamonds meeting $e$ join up cyclically to form a $2$--cell that intersects $e$ in its midpoint. For each face of $\mathcal{T}$, the two spinal half-edges meeting the barycenter of that face join to form a $1$--cell. The $0$--cells are simply the spinal vertices. Thus, there is a one-to-one correspondence between $i$--cells of $P$ and $(3-i)$--simplices of $\mathcal{T}$, for $i = 0, 1, 2$, with each corresponding cell and simplex pair intersecting transversely in their center/barycenter. 

Given $M\in\mathcal{M}_{g,k}$, Proposition \ref{prop-combinatorics} described the structure of its canonical triangulation $\mathcal{T}$. There is a single compact edge $e_0$, with incidence number $6g$. The corresponding $2$--cell of $P$ is made of $6g$ diamonds and has boundary a $6g$--gon. We will call this $2$--cell the \emph{big face} of $P$ and denote it by $G$. The remaining edges $e_1, \dotsc, e_k$ are non-compact and have incidence number $6$. Each exits a torus cusp (one per cusp) and the six diamonds meeting $e_i$ form a $2$--cell $H_i$, which we call a \emph{hexagonal face} of $P$. In $M$, the hexagonal face $H_i$ closes up into a torus representing the $i^{th}$ torus cusp. Recall that the torus cusp is made of two non-compact tetrahedra; each of them contributes three of the diamonds forming the hexagonal face. 


\section{Manifolds in \texorpdfstring{$\M_{k,k}$}{Mkk}}\label{sec:kkeven}


\begin{figure}
\centering
\includegraphics[width=1.75in,height=1.5in]{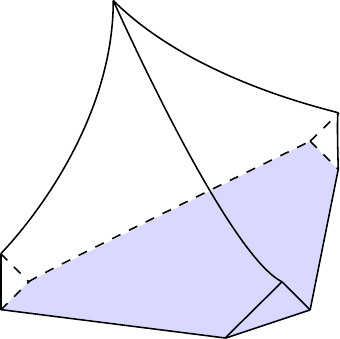}
\caption{The finite face of a partially truncated tetrahedron.}
\label{fig-finiteface}
\end{figure}


Given a partially truncated tetrahedron $\Delta^*$ with a unique ideal vertex, we call the unique hexagonal face which is not incident to the ideal vertex the \emph{finite face of $\Delta^*$} (see Figure \ref{fig-finiteface}).

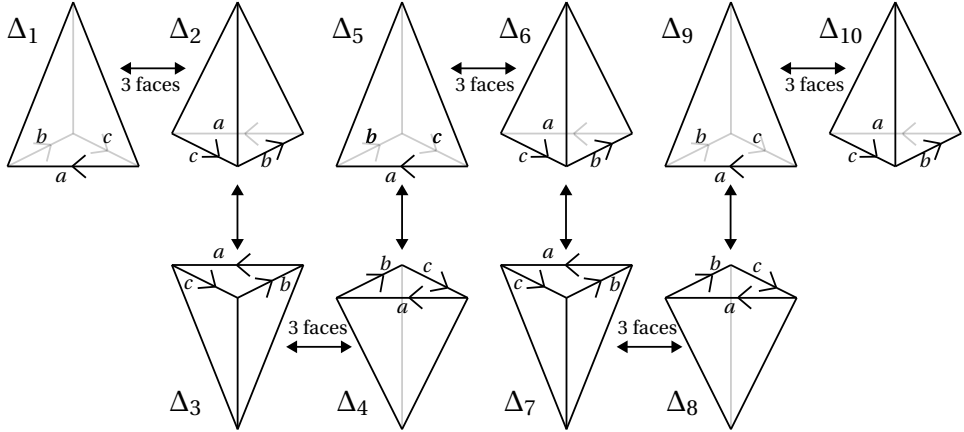
\begin{figure}[t]
\centering
\begin{overpic}[width=5in]{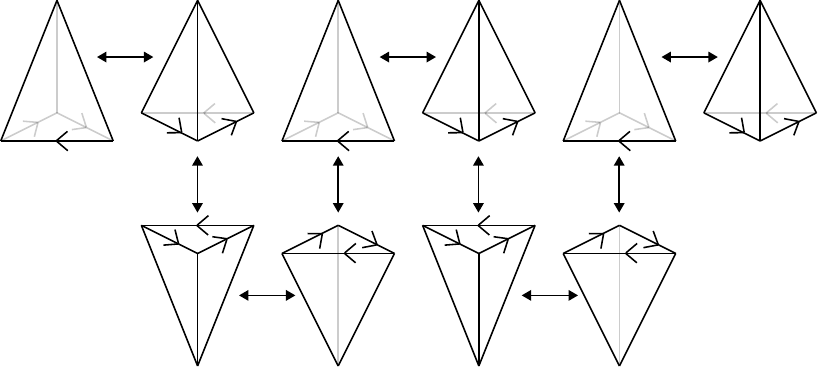}
\put(5, 26){\tiny$a$}
\put(39.5, 26){\tiny$a$}
\put(74, 26){\tiny$a$}
\put(21.5, 31.5){\tiny$a$}
\put(56, 31.5){\tiny$a$}
\put(90.5, 31.5){\tiny$a$}
\put(21.5, 18){\tiny$a$}
\put(56, 18){\tiny$a$}
\put(40.75, 12.25){\tiny$a$}
\put(75, 12.25){\tiny$a$}
\put(3, 30){\tiny$b$}
\put(10, 30){\tiny$c$}
\put(19, 28){\tiny$c$}
\put(18.5, 14.75){\tiny$c$}
\put(53, 14.75){\tiny$c$}
\put(28.5, 14.5){\tiny$b$}
\put(63, 14.5){\tiny$b$}
\put(26.5, 27.5){\tiny$b$}
 \put(53.5, 28){\tiny$c$}
 \put(61, 27.5){\tiny$b$}
 \put(88, 28){\tiny$c$}
 \put(95.5, 27.5){\tiny$b$}
\put(37.5, 30){\tiny$b$}
\put(44.5, 30){\tiny$c$}
\put(72, 30){\tiny$b$}
\put(78, 30){\tiny$c$}

\put(37.5, 30){\tiny$b$}
\put(44.5, 30){\tiny$c$}
\put(39, 16.75){\tiny$b$}
\put(43.5, 16.5){\tiny$c$}
\put(73.5, 16.75){\tiny$b$}
\put(78, 16.5){\tiny$c$}

\put(0, 41){$\Delta_1$}
\put(12,35.5){\tiny 3 faces}
\put(47,35.5){\tiny 3 faces}
\put(81.5,35.5){\tiny 3 faces}
\put(29.5,10){\tiny 3 faces}
\put(64,10){\tiny 3 faces}
\put(17, 41){$\Delta_2$}
\put(34, 41){$\Delta_5$}
\put(51.5, 41){$\Delta_6$}
\put(68.5, 41){$\Delta_9$}
\put(17, 2){$\Delta_3$}
\put(34.5, 2){$\Delta_4$}
\put(52, 2){$\Delta_7$}
\put(69, 2){$\Delta_8$}
\put(85, 41){$\Delta_{10}$}
\end{overpic}
\caption{A $5$-chain and its face identifications. The identifications respect the edge labels.}\label{fig-chain}
\end{figure}

\begin{definition}\label{def-chain}
An \emph{$\ell$-chain} is an ordered collection $\{\Delta_1,\Delta_2, \dotsc,\Delta_{2\ell}\}$ of $2\ell$ oriented non-compact tetrahedra, each of which has a single ideal vertex. Let $M(\Delta_1, \dotsc, \Delta_{2\ell})$ denote the $3$-manifold obtained from the partially truncated tetrahedra $\Delta_1^*, \dotsc, \Delta_{2\ell}^*$ by face pairings satisfying the following: 
\begin{itemize}
\item For $i=1,2,\dots,\ell$ each non-compact face of $\Delta_{2i-1}^*$ is glued to a non-compact face of $\Delta_{2i}^*$ such that all non-compact edges of  $\Delta_{2i-1}^*$ and  $\Delta_{2i}^*$ are identified to a single edge, and the union $\Delta_{2i-1}^* \cup \Delta_{2i}^*$ is homeomorphic to $\mathbb{T}^2 \times [0, \infty)$. 
\item For $i=1,2,\dots,\ell-1$ the finite face of $\Delta_{2i}^*$ is glued to the finite face of $\Delta_{2i+1}^*$
\end{itemize}
Moreover, relative to the induced orientations on faces of the $\Delta_i^*$, we require that the face-gluings are always by an orientation-reversing identification, so that $M(\Delta_1, \dotsc, \Delta_{2\ell})$
is oriented. See Figure \ref{fig-chain}. 
\end{definition}


\begin{remark}\label{remark-chain}
Recall that the faces of the $\Delta_i^*$ all have orientations. The finite face $\delta_1$ of $\Delta_1^*$ has edges with induced orientations; call them $a$, $b$, and $c$, in this order, along the oriented boundary of $\delta_1$ as in Figure \ref{fig-chain}. Under the identifications described above, the three edges of the finite face of $\Delta_i^*$ are identified with the three edges of the finite face of $\Delta_{i+1}^*$ for $i = 1, \dotsc, 2\ell -1$. Thus, all of these edges receive an orientation and a label $a$, $b$, or $c$ consistent with the gluings. 

The chain has two faces in its boundary: the finite faces $\delta_1$ of $\Delta_1$ and $\delta_{2\ell}$ of $\Delta_{2\ell}$. By definition, the cyclic ordering of the labels $a$, $b$, $c$ in the boundary of $\delta_1$ agrees with the induced orientation from $\delta_1$. Whether the same occurs in the oriented boundary of $\delta_{2\ell}$ depends on the parity of $\ell$; it agrees exactly when $\ell$ is odd. (This statement concerns only the cyclic ordering of labels, not the edge orientations.) 
\end{remark}


\begin{remark}\label{remark-edgecount}
Next we would like to glue $\delta_{2\ell}$ to $\delta_1$ by an orientation-reversing pairing to form an oriented $3$-manifold. Such a gluing will identify the finite edges of $M(\Delta_1, \dotsc, \Delta_{2\ell})$ by some bijection $\sigma\co \{a,b,c\} \to \{a,b,c\}$. By the previous remark, this bijection preserves the cyclic ordering of $a, b, c$ if and only if $\ell$ is even. 

For odd $\ell$, $\sigma$ fixes one label and transposes the other two. The resulting $3$-manifold will have two finite edges. For even $\ell$, $\sigma$ is either a $3$-cycle or the identity. The manifold will have either one or three finite edges, respectively. 
\end{remark}

\begin{figure}[t]
 \centering
\begin{tikzpicture}[scale=1.29]
\begin{scope}[decoration = {markings, mark = at position 0.65 with {\arrow{stealth}}}]
\draw[thick,lightgray] (.9,5.7) -- (.95,4.9); 
\draw[thick,lightgray,postaction={decorate}] (0.5,4.6) -- (.95,4.9);
\draw[thick,mycolor,postaction={decorate}] (.95,4.9) -- (1.3,4.6);
\draw[thick,postaction={decorate}] (1.3,4.6) -- (0.5,4.6);
\draw[thick] (.9,5.7) -- (0.5,4.6); 
\draw[thick] (.9,5.7) -- (1.3,4.6); 

\draw[thick,lightgray] (3.25,3.1) -- (3.3,4.3); 
\draw[thick,mycolor,postaction={decorate}] (2.85,4.0) -- (3.3,4.3);
\draw[thick,postaction={decorate}] (3.3,4.3) -- (3.65,4.0);
\draw[thick] (3.25,3.1) -- (2.85,4.0); 
\draw[thick] (3.25,3.1) -- (3.65,4.0); 
\draw[line width=3pt,white] (3.45,4.0) -- (3.05,4.0); 
\draw[thick,postaction={decorate}] (3.65,4.0) -- (2.85,4.0);

\draw[color=mycolor,thick,-{Stealth}] (3.1,4.5) .. controls (2.8,4.8) and (1.7,3.9) .. (1.2,4.3);
\node[color=mycolor] at (2.1,4.45) {\scriptsize $\sigma$};

\node at (.9,4.42) {\tiny $c$}; 
\node at (0.42,4.79) {\tiny $a$}; 
\node at (1.37,4.81) {\tiny $b$}; 

\node at (3.15,3.85) {\tiny $c$}; 
\node at (2.87,4.21) {\tiny $a$}; 
\node at (3.62,4.23) {\tiny $b$}; 

\node at (1.5,5.0) {$\dotsm$}; 
\node at (2.65,3.7) {$\dotsm$}; 
\node at (0.5,5.35) {\scriptsize $\Delta_1$};
\node at (3.75,3.4) {\scriptsize $\Delta_{2k}$};

\node[color=mycolor] at (0.77,3.55) 
{\tiny $\begin{aligned} \sigma \co &a \to b \\ 
&b \to c \\ 
&c \to a 
\end{aligned}$
};

\node at (2.25,5.55) {\footnotesize \underline{$M_L$}}; 
\draw[thick,lightgray] (0.2,2.9) -- (4.1,2.9) -- (4.1,5.9) --
(0.2,5.9) -- cycle;

\begin{scope}[shift={(5.0,0.0)}]
\draw[thick,lightgray] (.9,5.7) -- (.95,4.9); 
\draw[thick,lightgray,postaction={decorate}] (0.5,4.6) -- (.95,4.9);
\draw[thick,lightgray,postaction={decorate}] (.95,4.9) -- (1.3,4.6);
\draw[thick,mycolor,postaction={decorate}] (1.3,4.6) -- (0.5,4.6);
\draw[thick] (.9,5.7) -- (0.5,4.6); 
\draw[thick] (.9,5.7) -- (1.3,4.6); 

\draw[thick,lightgray] (3.25,3.1) -- (3.3,4.3); 
\draw[thick,mycolor,postaction={decorate}] (2.85,4.0) -- (3.3,4.3);
\draw[thick,postaction={decorate}] (3.3,4.3) -- (3.65,4.0);
\draw[thick] (3.25,3.1) -- (2.85,4.0); 
\draw[thick] (3.25,3.1) -- (3.65,4.0); 
\draw[line width=3pt,white] (3.45,4.0) -- (3.05,4.0); 
\draw[thick,postaction={decorate}] (3.65,4.0) -- (2.85,4.0);

\draw[color=mycolor,thick,-{Stealth}] (3.1,4.5) .. controls (2.8,4.8) and (1.7,3.9) .. (1.2,4.3);
\node[color=mycolor] at (2.1,4.5) {\scriptsize $\sigma^{-1}$};

\node at (.9,4.42) {\tiny $c$}; 
\node at (0.42,4.79) {\tiny $a$}; 
\node at (1.37,4.81) {\tiny $b$}; 

\node at (3.15,3.85) {\tiny $c$}; 
\node at (2.87,4.21) {\tiny $a$}; 
\node at (3.62,4.23) {\tiny $b$}; 

\node at (1.5,5.0) {$\dotsm$}; 
\node at (2.65,3.7) {$\dotsm$}; 
\node at (0.5,5.35) {\scriptsize $\Delta_1$};
\node at (3.75,3.4) {\scriptsize $\Delta_{2k}$};

\node[color=mycolor] at (0.77,3.55) 
{\tiny $\begin{aligned} \sigma^{-!} \co &a \to c \\ 
&b \to a \\ 
&c \to b 
\end{aligned}$
};

\node at (2.25,5.55) {\footnotesize \underline{$M_R$}}; 
\draw[thick,lightgray] (0.2,2.9) -- (4.1,2.9) -- (4.1,5.9) --
(0.2,5.9) -- cycle;
\end{scope}

\node at (7.15,2.64) {\rotatebox[origin=c]{-90}{$=$}};

\begin{scope}[shift={(5.0,-3.5)}]; 
\draw[thick] (.93,5.7) -- (.55,4.8); 
\draw[thick] (.93,5.7) -- (1.35,4.8); 
\draw[thick,lightgray,postaction={decorate}] (0.55,4.8) -- (1.35,4.8);
\draw[thick,postaction={decorate}] (0.9,4.5) -- (0.55,4.8);
\draw[thick,mycolor,postaction={decorate}] (1.35,4.8) -- (0.9,4.5);
\draw[line width=3pt,white] (.94,5.2) -- (0.9,4.6); 
\draw[thick] (.93,5.7) -- (0.9,4.5); 

\draw[thick] (3.27,3.1) -- (2.9,4.2); 
\draw[thick,postaction={decorate}] (2.9,4.2) -- (3.7,4.2);
\draw[thick,postaction={decorate}] (3.7,4.2) -- (3.25,3.9);
\draw[thick] (3.27,3.1) -- (3.7,4.2); 
\draw[thick] (3.27,3.1) -- (3.25,3.9); 
\draw[thick,mycolor,postaction={decorate}] (3.25,3.9) -- (2.9,4.2);

\draw[color=mycolor,thick,-{Stealth}] (3.1,4.5) .. controls (2.8,4.8) and (1.7,3.9) .. (1.2,4.3);
\node[color=mycolor] at (2.1,4.5) {\scriptsize $\sigma^{-1}$};

\node at (1.25,4.55) {\tiny $c$}; 
\node at (0.55,4.55) {\tiny $a$}; 
\node at (1.1,4.95) {\tiny $b$}; 

\node at (3.75,4.0) {\tiny $c$}; 
\node at (2.8,4.0) {\tiny $a$}; 
\node at (3.25,4.35) {\tiny $b$}; 

\node at (1.6,5.0) {$\dotsm$}; 
\node at (2.7,3.7) {$\dotsm$}; 
\node at (0.5,5.35) {\scriptsize $\Delta_1$};
\node at (3.75,3.4) {\scriptsize $\Delta_{2k}$};

\node[color=mycolor] at (0.77,3.55) 
{\tiny $\begin{aligned} \sigma^{-!} \co &a \to c \\ 
&b \to a \\ 
&c \to b 
\end{aligned}$
};

\node at (2.25,5.55) {\footnotesize \underline{$M_R$}}; 
\draw[thick,lightgray] (0.2,2.9) -- (4.1,2.9) -- (4.1,5.9) --
(0.2,5.9) -- cycle;
\end{scope}

\end{scope}
\end{tikzpicture}

 \caption{Constructing $M_L$ and $M_R$ from the chain $M(\Delta_1, \dotsc, \Delta_{2k})$ by gluing $\delta_{2k}$ to $\delta_1$. The gluing is completely determined by the highlighted pair of edges. }
 \label{fig-ML=MR}
 \end{figure}
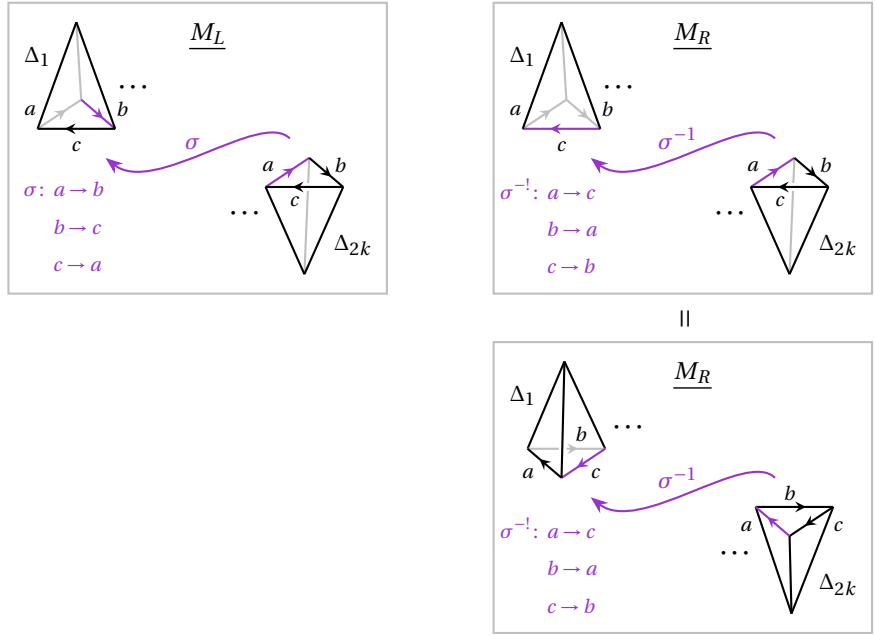

\begin{definition}\label{Twist}
Let $\ell$ be even. Let $\sigma$ be the $3$-cycle $(a \ b \ c)$. Define $M_L(\Delta_1, \dotsc, \Delta_{2\ell})$ to be obtained from $M(\Delta_1, \dotsc, \Delta_{2\ell})$ by gluing $\delta_{2\ell}$ to $\delta_1$ by a pairing inducing $\sigma$ on the edge labels. We refer to this gluing as the ``counterclockwise'' or ``left'' twist. Similarly, define $M_R(\Delta_1, \dotsc, \Delta_{2\ell})$ by using the ``clockwise'' or ``right'' twist $\sigma^{-1}$. See Figure \ref{fig-ML=MR}. Both $M_L$ and $M_R$ are oriented manifolds having a single finite edge. The quotient maps $M(\Delta_1, \dotsc, \Delta_{2\ell})\to M_L$ and $M(\Delta_1, \dotsc, \Delta_{2\ell}) \to M_R$ are denoted by $q_L$ and $q_R$, respectively. 
\end{definition}


\begin{proposition}\label{prop:chain}
Suppose $M \in \M_{k,k}$ and let $\Delta$ be any tetrahedron in $M$. Let $f \co \Delta_1 \to \Delta$ be any isomorphism of simplices taking ideal vertex to ideal vertex. Then there is a unique extension of $f$ to a combinatorial map 
\[ \bar{f} \co M(\Delta_1, \dotsc, \Delta_{2k}) \to M.\] 
Moreover, $k$ is even and $\bar{f}$ factors as 
\[ M(\Delta_1, \dotsc, \Delta_{2k}) \overset{q_{\ast}}{\longrightarrow} M_\ast(\Delta_1, \dotsc, \Delta_{2k}) \overset{h_f}{\longrightarrow} M\]
with $\ast = L$ or $R$ and $h_f$ a combinatorial isomorphism. 
\end{proposition}


\begin{proof}
$M$ admits an ideal triangulation $\mathcal{T}$ with $|\mathcal{T}|=2k$. 
By Proposition \ref{prop-combinatorics}(1), for each $i=1,\dots, k$, there are exactly two tetrahedra of $\mathcal{T}$ having $3$ vertices on $\Sigma_k$ and one vertex on the $i^{th}$ copy of $\mathbb{T}^2$. 
This exhausts all the tetrahedra in $\mathcal{T}$ since $|\mathcal{T}|=2k$ and we conclude that $\mathcal{T}$ contains no compact tetrahedra.
Note that each tetrahedron of $M$ shares a cusp with exactly one other tetrahedron, called its \emph{cusp neighbor}, and shares a finite face with exactly one other tetrahedron, called its \emph{finite face neighbor}. We define the map $\bar{f}$ inductively as follows. Supposing that $\bar{f}$ has been defined on $\Delta_1 \cup \dotsm \cup \Delta_i$, there are two cases. If $i$ is even then there is a unique map from $\Delta_{i+1}$ to the finite face neighbor of $\bar{f}(\Delta_i)$ extending $\bar{f}$ on the finite face of $\Delta_i$. If $i$ is odd then there is a unique extension sending $\Delta_{i+1}$ to the cusp neighbor of $\bar{f}(\Delta_i)$. 

This defines the map $\bar{f} \co M(\Delta_1, \dotsc, \Delta_{2k}) \to M$. Next we claim that each tetrahedron of $M$ is the image of exactly one $\Delta_i$, and hence $\bar{f}$ is surjective. It suffices to show that no two $\Delta_i$ have the same image in $M$. If we had $\bar{f}(\Delta_i) = \bar{f}(\Delta_j)$ for some $i<j$ then, since the image of $\Delta_r$ uniquely determines the image of $\Delta_{r+1}$, the images of $\Delta_i, \dotsc, \Delta_{j-1}$ would be the entire image of $\bar{f}$, and would comprise an open and closed submanifold of $M$, contradicting that $M$ is connected. Thus, $\bar{f}$ induces a bijection between the sets of tetrahedra, and is surjective. 

Finally, $\bar{f}$ is a quotient map and is injective on the complement of the finite faces $\delta_1$ and $\delta_{2k}$. The map $\bar{f}$ induces an identification between these faces, resulting in the triangulated manifold $M$. Since $M$ is orientable and has only one finite edge, the identification is one of the two possible twists creating $M_L$ or $M_R$, by Remark \ref{remark-edgecount}. 
\end{proof}


We are ready now to prove Theorem \ref{thm-kkeven}, which we restate for the readers convenience.  

\begin{figure}[t]
\centering
\includegraphics[scale=.1]{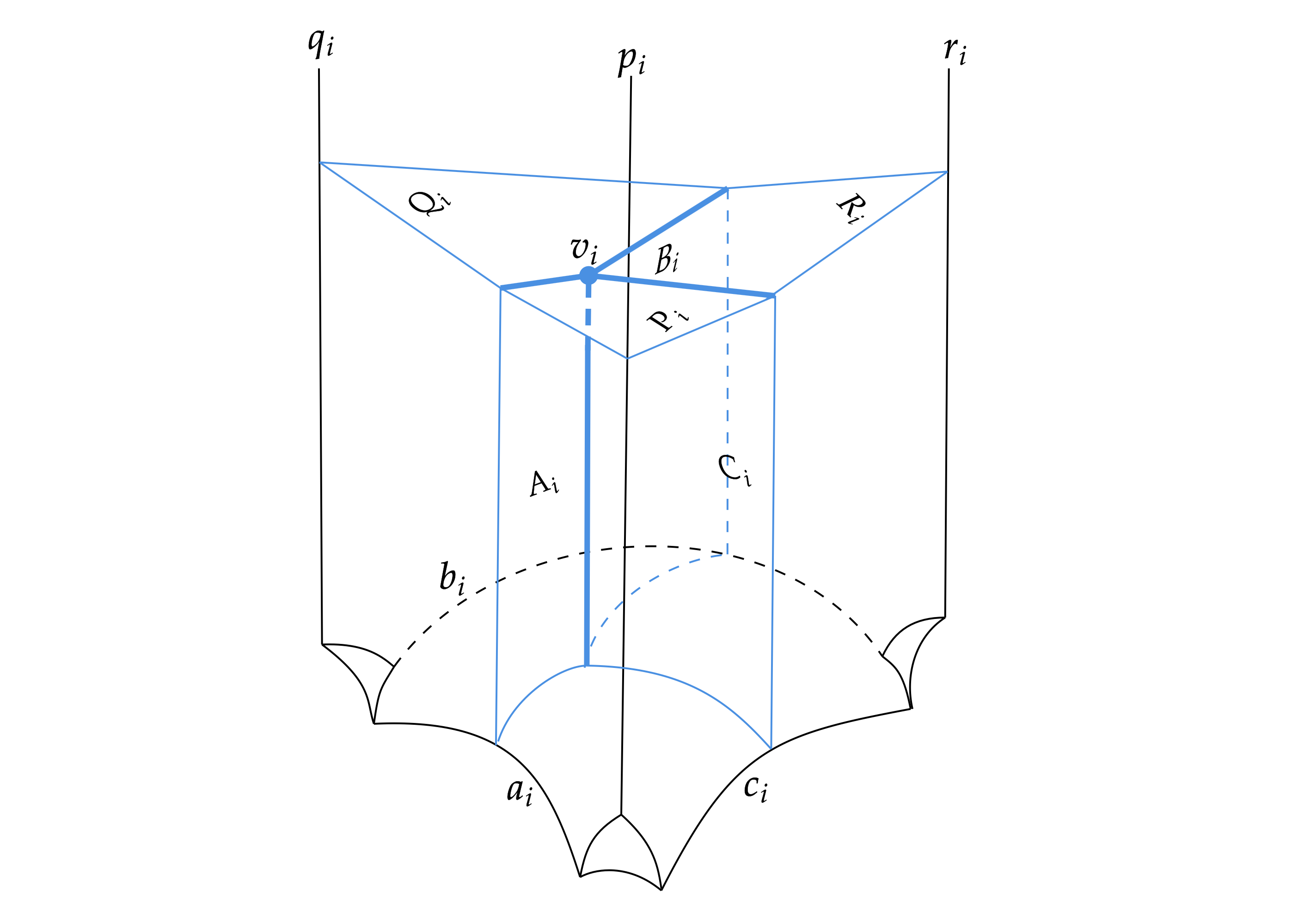}
\caption{Vertebra of a partially truncated tetrahedron with one ideal vertex.}
\label{fig-spine_notation}
\end{figure}


\thmkkeven*


\begin{proof} 
Proposition \ref{prop:chain} says that $\mathcal{M}_{k,k}$ consists of at most the homeomorphism types of $M_L$, $M_R$, so it suffices to show that $M_L$ and $M_R$ are indeed homeomorphic. 

Consider the constructions of $M_L$ and $M_R$ as shown in Figure \ref{fig-ML=MR}.  
The two faces $\delta_1, \delta_{2k}$ of the chain $M(\Delta_1, \dotsc, \Delta_{2k})$ are identified by a 1/3 rotation in the figure, in either direction. The second depiction of $M_R$ is combinatorially identical to the first; each tetrahedron in the chain has simply been rotated by 1/6 of a turn. There is a combinatorial isomorphism $M_L \to M_R$ taking each $\Delta_i$ to  $\Delta_i$ by an orientation-reversing isomorphism. This isomorphism may be seen by reflecting the left side of Figure \ref{fig-ML=MR} in the plane of the page (i.e. front-to-back), to obtain the second depiction of $M_R$. 
\end{proof}


At this point we digress to introduce some notation for the spine of the unique manifold in $\M_{k,k}$. This manifold admits two constructions, as $M_L$ and as $M_R$, and the spine notation we introduce depends on the choice of $M_L$ or $M_R$. 

The $k$-chain $\{\Delta_1, \Delta_2, \dots, \Delta_{2k}\}$ has three compact edge classes, which were labeled $a,b$, and $c$ (see Remark \ref{remark-chain}). We denote the compact edges of $\Delta_i$ by  $a_i, b_i$, and $c_i$ so that edges $a_{i}$ belong to the edge class $a$ and similarly for $b_i$, $c_i$. The non-compact edges of $\Delta_i$ are labeled $p_i, q_i$, and $r_i$ as shown in Figure \ref{fig-spine_notation}.
We denote the diamonds in the spine of $\Delta_i$ by $A_i, B_i, C_i, P_i, Q_i$, and $R_i$, where $A_i$ is the diamond meeting the edge $a_i$ and similarly for the others. We also denote the spinal vertex of $\Delta_i$ by $v_i$.

\begin{figure}[t]
\centering
\includegraphics[width=5in]{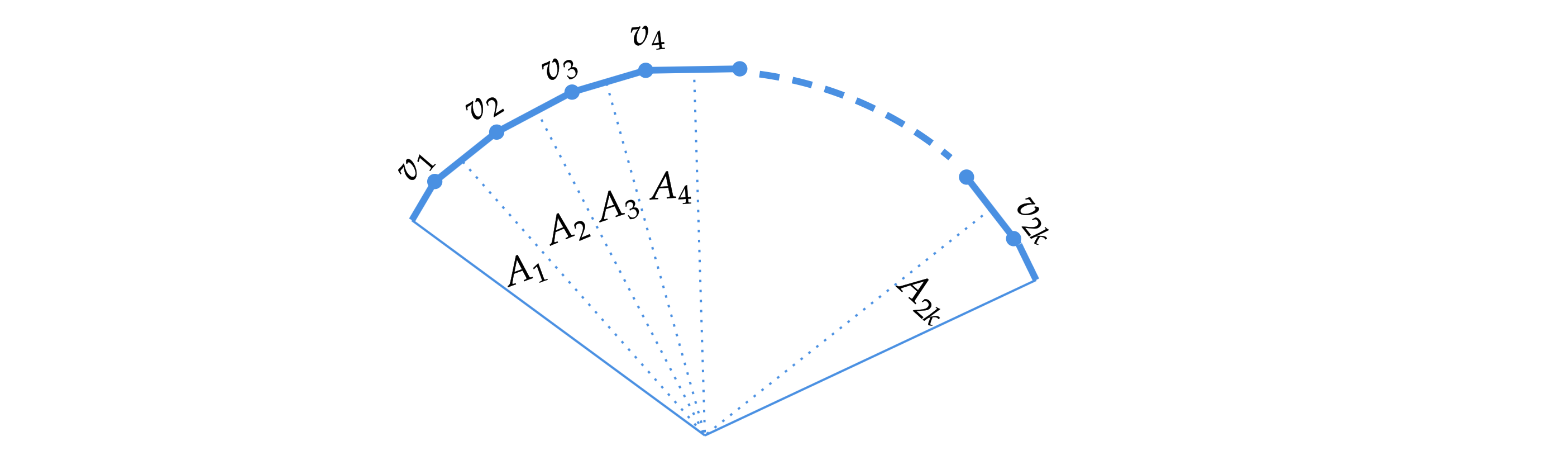}
\caption{Face pairings of a $k$-chain of partially truncated tetrahedra induce a chain of diamonds.}
\label{fig-chain of diamonds}
\end{figure}

\begin{figure}
\centering
\includegraphics[width=5.84in]{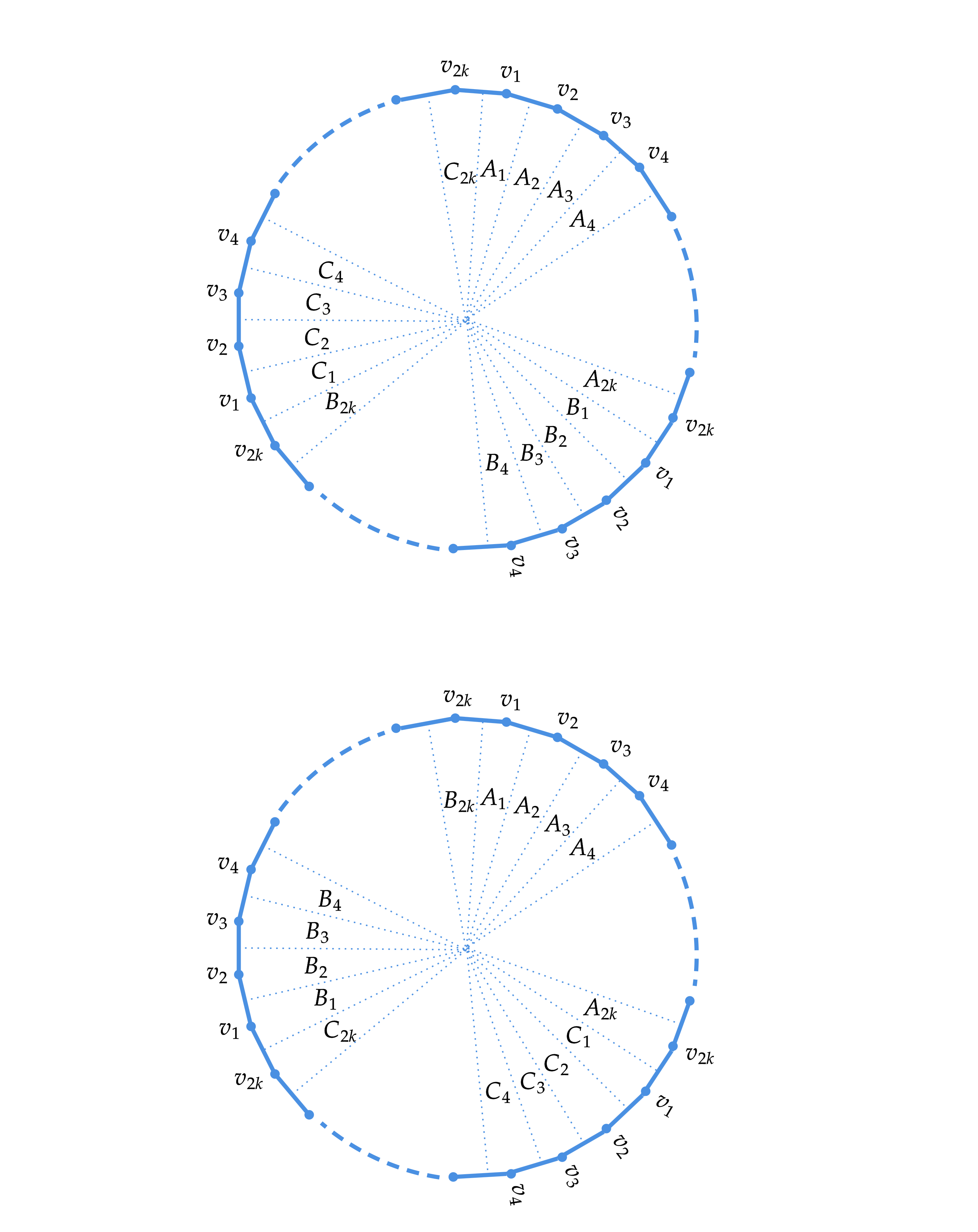}
\caption{The big face $G$ of the spine $P$ in $M_L$ (top) and in $M_R$ (bottom).}
\label{fig-big face G}
\end{figure}

The manifolds $M_L$ and $M_R$ were formed by closing up $M(\Delta_1, \dotsc, \Delta_{2k})$. In the construction of $M(\Delta_1, \dotsc, \Delta_{2k})$, all edges labeled $a_i$ were identified, all $b_i$ edges were identified, and all $c_i$ edges were identified. Thus, in $M(\Delta_1, \dotsc, \Delta_{2k})$, the diamonds $A_i$ meeting $a$ are all joined together, as are the $B_i$ diamonds and the $C_i$ diamonds. One may trace through the identifications shown in Figure \ref{fig-chain} to see that these diamonds are joined in the sequences $A_1, A_2, \dotsc, A_{2k}$; $B_1, \dotsc, B_{2k}$; and $C_1, \dotsc, C_{2k}$. See Figure \ref{fig-chain of diamonds} for the $A_i$ chain. The final face identification between $\delta_{2k}$ and $\delta_1$ creating $M_L$ or $M_R$ joins the three chains of diamonds into a single cycle, forming the big face $G$. In $M_L$ this cyclic ordering is
\[ A_1, \dotsc, A_{2k}, B_1, \dotsc, B_{2k}, C_1, \dotsc, C_{2k}\] 
and in $M_R$ the ordering is
\[ A_1, \dotsc, A_{2k}, C_1, \dotsc, C_{2k}, B_1, \dotsc, B_{2k}.\] 
See Figure \ref{fig-big face G}. 

With this terminology, we are now in a position to prove Theorem \ref{thm-isometries-kk}, which we restate for the reader's convenience.

\thmisometrieskk*

\begin{proof} 
Let $\mathcal{Z}$ be the set of isomorphisms of simplices $\Delta_1 \to \Delta_i$ (for $\Delta_i$ a simplex of $M_L$) taking ideal vertex to ideal vertex. Because $M_L \in \mathcal{M}_{k,k}$, Proposition \ref{prop:chain} says that each $f \in \mathcal{Z}$ yields a map $\bar{f}\co M(\Delta_1, \dotsc, \Delta_{2k}) \to M_L$ that factors as $h_f \circ q_L$ or $h_f \circ q_R$ with $h_f$ a combinatorial isomorphism $M_L \to M_L$ or $M_R \to M_L$. Define $\mathcal{Z}_L = \{f \in \mathcal{Z} \mid h_f \co M_L \to M_L \}$ and $\mathcal{Z}_R = \{f \in \mathcal{Z} \mid h_f \co M_R \to M_L \}$. 

Recall that the proof of Theorem \ref{thm-kkeven} exhibited a combinatorial isomorphism $M_L \to M_R$ taking $\Delta_1$ to $\Delta_1$ and reversing its orientation. Let $f_0$ be the inverse of that isomorphism restricted to $\Delta_1$. Then $f_0 \in \mathcal{Z}_R$, as $h_{f_0}$ is an isomorphism $M_R \to M_L$. Now pre-composition with $f_0$ gives a bijection $\mathcal{Z}_L \to \mathcal{Z}_R$, and therefore $\abs{\mathcal{Z}_L} = \abs{\mathcal{Z}_R} = \abs{\mathcal{Z}}/2 = 6k$. Finally note that $\Aut(M_L) = \{ h_f \mid f \in \mathcal{Z}_L\}$ and so $\abs{\Aut(M_L)} = 6k$. 

Next we turn to the spine of $M_L$. 
$\Aut(M_L)$ acts faithfully on the $2$--cell $G$ preserving the $6k$ vertices around its boundary \[v_1, \dotsc, v_{2k}, v_1, \dotsc, v_{2k}, v_1, \dotsc, v_{2k}.\] Also, since the action takes cusps to cusps and preserves the 
relation of being cusp neighbors, it takes each neighboring pair of vertices $(v_{2i-1}, v_{2i})$ to another such pair. The group of all such automorphisms of $G$ is the dihedral group $D_{3k}$. Thus our injective map $\Aut(M_L) \to \Aut(G)$ has image inside $D_{3k}$. Because $\Aut(M_L)$ and $D_{3k}$ have the same cardinality, we now have an isomorphism $\Aut(M_L) \cong D_{3k}$. 

It remains to determine which elements of $\Aut(M_L)$ are orientation-preserving. This is most easily understood via the action on the $2$-cell $G$. Recall that this cell is composed of the $6k$ diamonds $A_i, B_i, C_i$, arranged cyclically in the ordering
\[ A_1, \dotsc, A_{2k}, B_1, \dotsc, B_{2k}, C_1, \dotsc, C_{2k}.\] 
The group $\Aut(M_L)$, acting dihedrally on $G$, has a presentation 
\[ \langle \, r, t \mid r^2=1, \ t^{3k}=1, \ rtr=t^{-1} \, \rangle \]
with $t$ an order $3k$ rotation taking $A_1$ to $A_3$, $B_1$ to $B_3$, and $C_1$ to $C_3$; and $r$ a reflection exchanging $A_i$ with $A_{2k-i+1}$ and $B_i$ with $C_{2k-i+1}$ for all $i$. 

Now consider the isometries of $M_L$ extending $r$ and $t$ on $G$. The action of $t$ on $\Delta_1$ is determined by where $t$ takes the diamonds $A_1$, $B_1$, and $C_1$. Hence $\Delta_1$ is taken to $\Delta_3$, with edges labeled $a$, $b$, $c$ mapping to edges labeled $a$, $b$, $c$, respectively. Evidently this map is orientation-reversing (see Figure \ref{fig-chain}). Similarly, the isometry of $M_L$ extending $r$ takes $\Delta_1$ to $\Delta_{2k}$, with edges labeled $a$, $b$, $c$ mapping to edges labeled $a$, $c$, $b$, respectively. This map is orientation-preserving. The orientation-preserving subgroup $\Aut^+(M_L)$ has index 2 and contains $r$ and $t^2$, hence is equal to the dihedral subgroup $\langle \, r, t^2 \, \rangle \cong D_{3k/2}$. 
\end{proof}

\begin{remark}
The orientation-preserving isometries of $M_L$ do not always preserve the orientation of the $2$-cell $G$. Some elements (such as $r$) flip both the orientation of $G$ and its transverse orientation. 
\end{remark}


\section{Manifolds in \texorpdfstring{$\M_{k+1,k}$}{Mk+1k}}\label{sec:kplusone}

In this section, we prove Theorem \ref{thm-kplusone} by analyzing the combinatorial structure of the triangulations of manifolds in $\mathcal{M}_{k+1,k}$.
In particular, Proposition \ref{prop:even_chain_restriction} will prove parts (1) and (3) of Theorem \ref{thm-kplusone} and part (2) will be proved in Proposition \ref{prop:onehomeotype}.
The combination of these will complete the proof of Theorem \ref{thm-kplusone}.

By definition, any $M\in\mathcal{M}_{k+1,k}$ admits an ideal triangulation $\mathcal{T}$ of $2k+1$ tetrahedra. By Proposition \ref{prop-combinatorics}.(1), the following can be said about the tetrahedra in such a triangulation $\mathcal{T}$:
\begin{itemize}
\item  $\mathcal{T}$ has a compact tetrahedron $\Delta_0$ which has all four of its vertices on $\Sigma_{k+1}$.
\item For each $i=1,\dots,  k$ there are exactly two non-compact tetrahedra of $\mathcal{T}$ with $3$ vertices on $\Sigma_{k+1}$ and one vertex on the $i^{th}$ copy of $\mathbb{T}^2$. We denote these non-compact tetrahedra by $\Delta_{2i-1}$ and $\Delta_{2i}$.
\end{itemize}
Proposition \ref{prop-combinatorics}.(2) prescribes that the three non-compact faces of $\Delta_{2i-1}$ must be identified with the three non-compact faces of $\Delta_{2i}$.
All the remaining faces in $\mathcal{T}$ are compact, four from $\Delta_0$ and one from each $\Delta_i$ with $i\ge 1$.
The following lemma describes the necessary identifications of the remaining compact faces in $\mathcal{T}$ in order to construct a manifold in $\M_{k+1,k}$.
In what follows, as well as the remainder of the subsection, we will also abusively refer to the manifold $M(\Delta_1,\dots,\Delta_{2\ell})$ from Definition \ref{def-chain} as an $\ell$-chain.


\begin{lemma}\label{lemma-chains of M_k+1,k}
Let $M\in\mathcal{M}_{k+1,k}$ and $\mathcal{T}$ with notation as above, then $M$ must be obtained in one of the following two ways:
\begin{enumerate}
\item Gluing each finite face of a single $k$-chain, $\{\Delta_1,\dots,\Delta_{2k}\}$, to two faces of $\Delta_0$ and identifying the remaining two compact faces of $\Delta_0$ with each other.
\item Gluing each finite face of $\Delta_0$ to the finite faces of two distinct chains, $\{\Delta_1,\dots,\Delta_{2\ell}\}$, $\{\Delta_{2\ell+1},\dots,\Delta_{2k}\}$, for some $1\le \ell<k$.
\end{enumerate}
\end{lemma}


\begin{proof}
We prove both cases simultaneously.
As above, each pair $(\Delta_{2i-1},\Delta_{2i})$ is glued on non-compact faces and has two remaining compact faces.
One can encode the connectedness of $M$ in the adjacency graph of these pairs.
Specifically, let $\mathcal{G}_M$ be the finite graph with $k+1$ vertices $x_0,\dots, x_{k}$, with one vertex $x_0$ for $\Delta_0$, one vertex $x_i$ per pair $(\Delta_{2i-1},\Delta_{2i})$, and edges denoting identification of compact faces.
For each face of $\Delta_0$ which is identified with another face of $\Delta_0$, we give the vertex $x_0$ in $\mathcal{G}_M$ a self-loop, and similarly if the two compact faces of a pair $(\Delta_{2i-1},\Delta_{2i})$ are identified, then we give $x_i$ a self-loop.
By construction, it follows that $\mathcal{G}_M$ has the property that $x_0$ has degree $4$ and each $x_i$ has degree $2$ for $i\ge 1$ (our convention is that a self-loop adds $2$ to local degree).

From this description one can readily see that $\mathcal{G}_M$ is connected if and only if it is a wedge of two cycles with wedge point $x_0$ (here we count a self-loop at $x_0$ as a cycle).
Indeed, choose any edge incident to $x_0$ which connects to some vertex $x_i\neq x_0$ and follow the resulting path.
Necessarily one such exists unless $x_0$ is the only vertex in its connected component. 
If $\mathcal{G}_M$ is connected, every vertex $x_j$ that one encounters along this path is degree $2$ with no self-loops until one returns to $x_0$.
This must happen by returning on a different edge incident to $x_0$ as otherwise there was a vertex of degree $\ge3$ at some point along the path.
At this point, either one has traversed all edges of the graph except a single self-loop based at $x_0$ or one may repeat this procedure using one of the remaining two edges incident to $x_0$ and conclude the similar result for that path. 
In the former case, relabeling the tetrahedra to match the order of this cycle, we see that $M$ is constructed by identifying the finite faces of the chain $\{\Delta_1,\dots,\Delta_{2k}\}$ with those of $\Delta_0$ as well as two finite faces of $\Delta_0$ with themselves.
In the latter case, relabeling the tetrahedra to match the order of these two cycles, we conclude that $M$ is constructed by identifying the finite faces of the chains $\{\Delta_1,\dots,\Delta_{2\ell}\}$, $\{\Delta_{2\ell+1},\dots,\Delta_{2k}\}$ of lengths $\ell$, $k-\ell$ (respectively) for some $1\le \ell<k$.
This completes the proof.
\end{proof}


\begin{corollary}\label{cor-chains of M_k+1,k}
For any $M \in \M_{k+1,k}$ there is a canonical partition of the set of torus components of $\partial M$ into two sets (one of which may be empty) corresponding to the cycles from Cases (1) and (2) of Lemma \ref{lemma-chains of M_k+1,k}.
Therefore the cardinalities of these partitions (equivalently cycles in $\mathcal{G}_M$), $i(M)$, $j(M)$, are canonical invariants associated to any $M\in\mathcal{M}_{k+1,k}$.
We always order these invariants so that $0\le i(M)\leq j(M)$ and $i(M)=0$ if and only if we are in Case (1).
\end{corollary}


As we will see below, Corollary \ref{cor-chains of M_k+1,k} is necessary and almost sufficient to construct a manifold $M\in\M_{k+1,k}$, however, we have not yet taken into account the restrictions arising from Proposition \ref{prop-combinatorics}.
Importantly, Proposition \ref{prop-combinatorics}.(2) dictates that, after identification, there is precisely one edge class incident to the boundary $\Sigma_{k+1}$.
That is to say, all finite edges on the finite faces have been identified with each other.
The content of the next lemma is that such a construction is always possible unless $i(M),j(M)$ from Corollary \ref{cor-chains of M_k+1,k} are both even.


\begin{proposition}\label{prop:even_chain_restriction}
Let $M\in\mathcal{M}_{k+1,k}$. Then at least one of $i(M)$, $j(M)$ from Corollary \ref{cor-chains of M_k+1,k} is odd.
Additionally, for every pair of integers $(i, j)$ such that $0\le i\le j\le k$, $i+j=k$, and at least one of $i$, $j$ is odd, there is an associated manifold $M\in\mathcal{M}_{k+1,k}$ such that $(i,j)=(i(M),j(M))$.
\end{proposition}


\begin{proof}
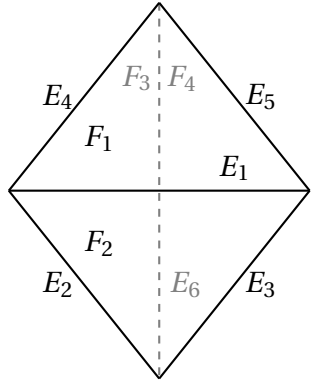
\begin{figure}[t]
\begin{tikzpicture}
    \coordinate (A) at (0,0);
    \coordinate (B) at (4,0);
    \coordinate (C) at (2,2.5);
    \coordinate (D) at (2,-2.5);
    
    \draw[thick] (A) -- (B) node[pos=.75, above] {$E_1$};
    \draw[thick] (A) -- (C) node[midway, left] {$E_4$};
    \draw[thick] (B) -- (C) node[midway, right] {$E_5$};
    \draw[thick] (A) -- (D) node[midway, left] {$E_2$};
    \draw[thick] (B) -- (D) node[midway, right] {$E_3$};
    \draw[dashed,thick,opacity=.5] (C) -- (D) node[pos=.75, right] {$E_6$};
    
    \node at (1.2,0.7) {$F_1$};
    \node at (1.2,-0.7) {$F_2$};
    \node[opacity=.5] at (1.7,1.5) {$F_3$};
    \node[opacity=.5] at (2.3,1.5) {$F_4$};
\end{tikzpicture}
\caption{A labeling of the faces and edges of the compact tetrahedron, $\Delta_0$. The faces $F_3$, $F_4$ as well as the edge $E_6$ are on the back side of this tetrahedron.}\label{fig:compactface}
\end{figure}
We will prove both statements simultaneously. 
We denote the chains of lengths $i(M)$, $j(M)$ (respectively) furnished by Lemma \ref{lemma-chains of M_k+1,k} as
$$\mathcal{K}_1=\{\Delta_1,\dots,\Delta_{2i(M)}\}, \quad \mathcal{K}_2=\{\Delta_{2i(M)+1},\dots,\Delta_{2k}\},$$
which are glued to the compact tetrahedron $\Delta_0$.
We allow for the possibility that $i(M)=0$, in which case the former chain is the empty chain which still identifies two faces of $\Delta_0$ as in Lemma \ref{lemma-chains of M_k+1,k}(1).

Enumerate the edges and faces of $\Delta_0$ as in Figure \ref{fig:compactface}, so that $F_1$, $F_2$ share the edge $E_1$ and $F_3$, $F_4$ share the edge $E_6$.
Without loss of generality, we may assume that the unique compact face of $\Delta_1$, $\Delta_{2i(M)}$, $\Delta_{2i(M)+1}$, $\Delta_{2k}$ is identified with $F_1$, $F_2$, $F_3$, $F_4$, respectively.
Moreover, Proposition \ref{prop-combinatorics} shows that $M\in\mathcal{M}_{k+1,k}$ if and only if after identification there is a single edge class shared among the compact faces of $\Delta_0,\dots,\Delta_{2k}$.
This is equivalent to asking that, after identification, the edges $E_1,\dots, E_6$ of $\Delta_0$ lie in a single edge class.
In the remainder of the proof, we will analyze identifications coming from chains of even and odd lengths and how they affect edge classes.

We first make some general reductions.
Note that the identification of $F_1$ with the compact face of $\Delta_1$, propagates the edge class of $E_1$, $E_2$, $E_3$ throughout the compact faces of the chain $\mathcal{K}_1$, ultimately endowing the compact face of $\Delta_{2i(M)}$ with a labeling of these edge classes.
Upon identification of $\Delta_{2i(M)}$ with $F_2$, the edge classes of $E_1$, $E_2$, $E_3$ then get identified with those of $E_1$, $E_4$, $E_5$.
There are six possibilities for this identification, indexed by elements of $S_3$, three of which are possible when the chain has even length and three of which are possible when the chain has odd length, see Remark \ref{remark-edgecount}.
We note that two of these identifications can be ruled out immediately, namely those identifications which take $E_1$ to $E_1$.
Indeed if this was the case, then the edge class of $E_1$ is  given by $[E_1]=\{E_1\}$.
Since the remaining identifications of edges of $F_3$ with edges of $F_4$ via $\mathcal{K}_2$ only identify edges contained in the set $\{E_2,\dots,E_6\}$, there will be more than one edge class resulting from identifications via $\mathcal{K}_1$, $\mathcal{K}_2$, contradicting Proposition \ref{prop-combinatorics}.
In the remaining cases, $E_1$ always gets identified with another edge, therefore we reduce to studying the edge classes $[E_2],\dots, [E_6]$, as $E_1$ is necessarily contained in one of these classes.
The identical reductions hold for the edge $E_6$, the faces $F_3$, $F_4$, and their identifications via $\mathcal{K}_2$.
Consequently, we reduce to studying the edge classes of the four edges $E_2$, $E_3$, $E_4$, $E_5$ provided we only consider the identifications via $\mathcal{K}_1$ (resp. $\mathcal{K}_2$) which do not identify $E_1$ (resp. $E_6$) with itself.


\begin{figure}[t]
    \centering
    \subfigure{
        \begin{tikzpicture}[scale=1.5]
            \draw[red, thick] (0,1) -- (0,0);
            \draw[red, thick] (1,1) -- (1,0);
            \filldraw[black] (0,1) circle (2pt) node[anchor=east] {$E_4$};
            \filldraw[black] (1,1) circle (2pt) node[anchor=west] {$E_5$};
            \filldraw[black] (0,0) circle (2pt) node[anchor=east] {$E_2$};
            \filldraw[black] (1,0) circle (2pt) node[anchor=west] {$E_3$};
        \end{tikzpicture}
    }\hspace{.5in}
    \subfigure{
        \begin{tikzpicture}[scale=1.5]
            \draw[red, thick] (0,1) -- (1,0);
            \draw[red, thick] (1,1) -- (0,0);
            \filldraw[black] (0,1) circle (2pt) node[anchor=east] {$E_4$};
            \filldraw[black] (1,1) circle (2pt) node[anchor=west] {$E_5$};
            \filldraw[black] (0,0) circle (2pt) node[anchor=east] {$E_2$};
            \filldraw[black] (1,0) circle (2pt) node[anchor=west] {$E_3$};
        \end{tikzpicture}
    }\hspace{.5in}
    \subfigure{
        \begin{tikzpicture}[scale=1.5]
            \draw[red, thick] (0,1) -- (1,1);
            \draw[red, thick] (0,0) -- (1,0);
            \filldraw[black] (0,1) circle (2pt) node[anchor=east] {$E_4$};
            \filldraw[black] (1,1) circle (2pt) node[anchor=west] {$E_5$};
            \filldraw[black] (0,0) circle (2pt) node[anchor=east] {$E_2$};
            \filldraw[black] (1,0) circle (2pt) node[anchor=west] {$E_3$};
        \end{tikzpicture}
    }
    \caption{The edge identifications induced by chains of even and odd length. The left figure is the edge identification induced by the chain $\mathcal{K}_1$ if it has odd length, the middle figure is the edge identification induced by any chain of even length, and the right figure is the edge identification induced by$\mathcal{K}_2$ if it has odd length.}\label{fig:edgeclass}
\end{figure}
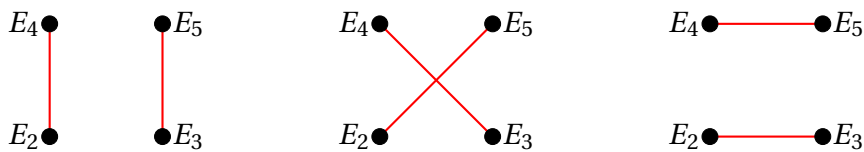


The key observation is then the following: if $\mathcal{K}_i$ is of even length, then the identification of the two compact faces $F_{2i-1},F_{2i}$ of $\Delta_0$ that $\mathcal{K}_i$ is incident to identifies edges diagonally across from each other, whereas if $\mathcal{K}_i$ is of odd length, then it identifies the edges which are adjacent and share a common vertex with the central edge, $E_{5i-4}$.
Indeed, this follows immediately from analyzing how edge labels propagate under even or odd chains (see Figures \ref{fig-chain} and \ref{fig-ML=MR}).
This is perhaps most easily visualized in the adjacency graph in Figure \ref{fig:edgeclass}, where vertices correspond to the edges $E_2$, $E_3$, $E_4$, $E_5$ and edges correspond to identification via a chain.
We note that $\mathcal{K}_1$, $\mathcal{K}_2$ being even length induce the same identification on the edges, however, this is not true for odd length chains.
The key difference is that the central edge, $E_1$, shared by $F_1$, $F_2$ is in a different orientation as the central edge, $E_6$, shared by $F_3$, $F_4$, which explains the discrepancy in these two cases in Figure \ref{fig:edgeclass}.

Performing the edge class identifications from $\mathcal{K}_1$, $\mathcal{K}_2$ corresponds to taking one of the left or center graphs from Figure \ref{fig:edgeclass} (which are the identifications via $\mathcal{K}_1$) and layering it on top of one of the center or right graphs from Figure \ref{fig:edgeclass} (which are the identifications via $\mathcal{K}_2$).
The resulting graph is connected if and only if there is a single edge class after all identifications have been performed.
From this one can immediately see that this happens if and only if at least one of the chains is odd length (equivalently, at least one of the right or left graphs is used).
This concludes the proof.
\end{proof}


At this point, we have completed the proofs of Parts (1) and (3) of Theorem \ref{thm-kplusone}.
Moreover, the proof of Proposition \ref{prop:even_chain_restriction} yields the following, which will be useful for proving Part (2).


\begin{corollary}\label{cor:atmostfour}
Fix a pair $(i, j)$ such that $0\le i\le j\le k$, $i+j=k$, and at least one of $i$, $j$ is odd.
Then there are at most four homeomorphism types of $M\in\mathcal{M}_{k+1,k}$ such that $i=i(M)$ and $j=j(M)$.
Moreover, all such manifolds are built from two chains $\mathcal{K}_1$, $\mathcal{K}_2$ as in Lemma \ref{lemma-chains of M_k+1,k} (where one is possibly empty) and each possible homeomorphism type arises from choosing one of two possible identifications for a single finite face of $\mathcal{K}_i$.
\end{corollary}


\begin{proof}
This follows from the confluence of Corollary \ref{cor:combinatorialiso} and the proof of Proposition \ref{prop:even_chain_restriction}.
Indeed, any two such manifolds are homeomorphic if and only if they are combinatorially isomorphic. 
The proof of Proposition \ref{prop:even_chain_restriction} then shows that there are at most four possible combinatorial types, corresponding to the two admissible face identifications via $\mathcal{K}_1$ and similarly the two admissible face identifications via $\mathcal{K}_2$.
\end{proof}


To complete the proof of Theorem \ref{thm-kplusone}, we show that there is in fact only one homeomorphism type.


\begin{proposition}\label{prop:onehomeotype}
Suppose that $M,M'\in\mathcal{M}_{k+1,k}$ such that $i(M)=i(M')$ and $j(M)=j(M')$, then $M$ is homeomorphic to $M'$.
In particular, the invariants $i(M), j(M)$ are complete invariants of the homeomorphism type of $M$.
\end{proposition}


\begin{proof}
By Corollary \ref{cor:atmostfour}, for any fixed pair $(i,j)$, there are at most four elements in $\mathcal{M}_{k+1,k}$ with these invariants, which we enumerate as $M_1,\dots, M_4$.
We will show that there are combinatorial isomorphisms between each pair $(M_i, M_j)$ which will yield the desired result.

Note that if $M_i$, $M_j$ are combinatorially isomorphic then such an isomorphism must preserve $\Delta_0$.
Moreover, because this isomorphism must take chains to chains, it will either fix $F_1\cup F_2$ and $F_3\cup F_4$ (setwise) or interchange the two.
When $i(M_i)\neq j(M_i)$, this isomorphism cannot do the latter and therefore must fix the edges $E_1$, $E_6$ (setwise). 
When $i(M_i)=j(M_i)$, it is possible to have automorphisms that interchange $E_1$ and $E_6$, however, we will only require the more restrictive class of combinatorial isomorphisms that fix these $E_1$, $E_6$ in order to prove the result.

Let $X_0$ denote the set of combinatorial automorphisms of $\Delta_0$ preserving both $E_1$, $E_6$ (as sets), which is isomorphic to the Klein four group $\Z/2\Z\times \Z/2\Z$.
The group $X_0$ is generated by a reflection in the hyperplane containing $E_1$ and the midpoint of $E_6$, denoted $R_1$, and a reflection in the hyperplane containing $E_6$ and the midpoint of $E_1$, denoted $R_6$. 
One easily checks that $R_1R_6=R_6R_1$ is a rotation of $\Delta_0$ by $\pi$ and interchanges the position of $F_1$ with $F_2$ and $F_3$ with $F_4$.

Let $\mathcal{M}=\{M_1,M_2,M_3,M_4\}$ and fix some $M_i\in\mathcal{M}$.
Then we claim that any element $h\in X_0$ can be extended to a combinatorial isomorphism $f_h:M_i\to M_j$ for some $M_j\in\mathcal{M}$ with the property that $f_h\vert_{\Delta_0}= h\cdot\Delta_0$, where $j=j(h)$ depends on $h$.
Indeed, it suffices to show this for $R_1$, $R_6$ provided that $f_{R_1}\circ f_{R_6}=f_{R_6}\circ f_{R_1}$, as then $h=R_1R_6$ follows by composition.
To this end, we define the following two operations on chains:
\begin{enumerate}
\item[] \emph{Inversion:} Given a chain $\mathcal{K}=\{\Delta_1,\dots,\Delta_{2\ell}\}$ of length $\ell$, we let $\mathcal{K}^I=\{\Delta_{2\ell},\dots,\Delta_{1}\}$ which is the chain with reversed ordering on the tetrahedra.\smallskip
\item[] \emph{Face Reflection:} Given a tetrahedron $\Delta_i$, we define $\overline{\Delta}_i$ to be the tetrahedron obtained by reflecting in the unique hyperplane through the ideal vertex, the original edge adjacent to $a$, $b$, and bisecting the edge $c$ (with labels as in Figure \ref{fig-chain}).
Given a chain $\mathcal{K}=\{\Delta_1,\dots,\Delta_{2\ell}\}$ of length $\ell$, we define the face reflected chain by $\overline{\mathcal{K}}=\{\overline{\Delta}_1,\dots,\overline{\Delta}_{2\ell}\}$. Note that reflecting in all faces is compatible with gluing so $\overline{\mathcal{K}}$ is also a chain of length $\ell$.
\end{enumerate}
One checks that the operations of inversion and face reflection commute, so we may unambiguously write $\overline{\mathcal{K}}^I$.
Enumerating the chains of $M_i$ as $\mathcal{K}_1$, $\mathcal{K}_2$ as before and writing $M(\mathcal{K})$ to denote the correspond manifold as in Definition \ref{def-chain}, we define the following maps
\begin{align*}
f_{R_1}:M_i&\to M_{j(R_1)},\\
\Delta_0&\mapsto R_1\cdot \Delta_0,\\
M(\mathcal{K}_1)&\mapsto M(\overline{\mathcal{K}}^I_1),\\
M(\mathcal{K}_2)&\mapsto M(\overline{\mathcal{K}}_2),
\end{align*}
and
\begin{align*}
f_{R_6}:M_i&\to M_{j(R_6)}\\
\Delta_0&\mapsto R_6\cdot \Delta_0,\\
M(\mathcal{K}_1)&\mapsto M(\overline{\mathcal{K}}_1),\\
M(\mathcal{K}_2)&\mapsto M(\overline{\mathcal{K}}^I_2).
\end{align*}
These maps are continuous by the pasting lemma, as the restrictions to the faces of $\Delta_0$ agree after gluing.
It is straightforward to check that they are their own inverses and hence combinatorial involutions. 
It is also straightforward to check that $f_{R_1}\circ f_{R_6}=f_{R_6}\circ f_{R_1}$, as $R_1 R_6\cdot\Delta_0=R_1\cdot( R_6\cdot\Delta_0)=R_6\cdot (R_1\cdot\Delta_0)$ and, in either order of composition, the induced maps on chains are given by $M(\mathcal{K}_i)\mapsto M(\mathcal{K}^I_i)$.
Moreover, as $M_j$ is built by gluing chains of the same length as $M_i$, necessarily $M_j$ is an element of $\mathcal{M}$, as required.

Finally, to conclude the result it suffices to prove that the action of $X_0$ via $f_h$ on the set $\mathcal{M}$ is transitive.
Equivalently, we show that the action on this set is free; that is, if $f_h(M_i)=M_{j(h)}=M_i$ then $h=\mathrm{id}_{\Delta_0}$.

Suppose that $f_h(M_i)=M_i$ for some $i=1,\dots, 4$.
As $f_h$ is a combinatorial isomorphism it must induce a homeomorphism of $\Delta_0$ with itself and therefore induces a homeomorphism of the big face $G$ of the spine, as described in Section \ref{subsection:spinebackground}.
Indeed in this setting, the spinal vertex of $\Delta_0$ corresponds to the central vertex of $G$ and hence $f_h$ induces a homeomorphism of all of the spinal edges incident to this vertex and consequently all of the diamonds incident to this spinal vertex as well.
Moreover, there are six diamonds internal to $\Delta_0$ and six chains of diamonds arising from the chains $\mathcal{K}_1$, $\mathcal{K}_2$; three chains of diamonds for each of $\mathcal{K}_1$, $\mathcal{K}_2$, which are the similar chains to those described in Section \ref{sec:kkeven}.
The difference here is that the chains of diamonds are glued on non-spinal edges to a diamond of $\Delta_0$ via the face identifications coming from the $\mathcal{K}_i$ and therefore the big face alternates between a diamond of $\Delta_0$ and a chain of diamonds from $\mathcal{K}_1$, $\mathcal{K}_2$.

Figure \ref{fig:diamonds} gives a description of the big face $G$.
In this figure, the grey slices correspond to diamonds of $\Delta_0$, the labeled slices are chains of diamonds with subscript corresponding to the pair of faces they connect, superscript enumerating which of the three chains of diamonds it is, and arrows corresponding to whether the chain of diamonds is increasing -- connecting a face $F_r$ with $F_{r+1}$ -- or decreasing -- connecting a face $F_{r+1}$ with $F_{r}$.
Drawing such a diagram involves choosing both a starting chain, that is, the top rightmost label, as well as a direction to travel through this chain.
These choices correspond, respectively, to a cyclic rotation of the diagram and a reflection in the plane and therefore this diagram is an invariant of $M_i$ up to these two diagrammatic automorphisms.

\begin{figure}[t]
\begin{tikzpicture}[>=stealth]
  \def\n{12}              
  \def\radius{3cm}          
  \def\slice{360/\n}       
  \def\arcr{3.2cm}          
  \def\archo{10}             

  \foreach \i in {0,...,11} {
    \pgfmathsetmacro{\ang}{\i*\slice}
    \pgfmathtruncatemacro{\m}{mod(\i,2)}
    \ifnum\m=0
      \def\shadecolor{gray!50}
    \else
      \def\shadecolor{white}
    \fi
    \draw[fill=\shadecolor,draw=black] (0,0) -- (\ang:\radius) arc (\ang:\ang+\slice:\radius) -- cycle;
  }

  \pgfmathsetmacro{\angmid}{1*\slice + \slice/2}
  \node at (\angmid:\radius*0.6) {$F_{1/2}^{(1)}$};
  \draw[->] (\angmid+\archo:\arcr) arc[start angle=\angmid+\archo, end angle=\angmid-\archo, radius=\arcr];

  \pgfmathsetmacro{\angmid}{11*\slice + \slice/2}
  \node at (\angmid:\radius*0.6) {$F_{3/4}^{(1)}$};
  \draw[->] (\angmid-\archo:\arcr) arc[start angle=\angmid-\archo, end angle=\angmid+\archo, radius=\arcr];

  \pgfmathsetmacro{\angmid}{9*\slice + \slice/2}
  \node at (\angmid:\radius*0.6) {$F^{(2)}_{1/2}$};
  \draw[->] (\angmid-\archo:\arcr) arc[start angle=\angmid-\archo, end angle=\angmid+\archo, radius=\arcr];

  \pgfmathsetmacro{\angmid}{7*\slice + \slice/2}
  \node at (\angmid:\radius*0.6) {$F^{(3)}_{1/2}$};
  \draw[->] (\angmid-\archo:\arcr) arc[start angle=\angmid-\archo, end angle=\angmid+\archo, radius=\arcr];

  \pgfmathsetmacro{\angmid}{5*\slice + \slice/2}
  \node at (\angmid:\radius*0.6) {$F^{(2)}_{3/4}$};
  \draw[->] (\angmid-\archo:\arcr) arc[start angle=\angmid-\archo, end angle=\angmid+\archo, radius=\arcr];

  \pgfmathsetmacro{\angmid}{3*\slice + \slice/2}
  \node at (\angmid:\radius*0.6) {$F^{(3)}_{3/4}$};
  \draw[->] (\angmid+\archo:\arcr) arc[start angle=\angmid+\archo, end angle=\angmid-\archo, radius=\arcr];

\end{tikzpicture}
\caption{One of the four possible configurations of $F_{1/2}^{(1)}F_{3/4}^{(1)}F^{(2)}_{1/2}F^{(3)}_{1/2}F^{(2)}_{3/4}F^{(3)}_{3/4}$. In this configuration, $F_{1/2}$ is increasing and $F_{3/4}$ is decreasing. 
Each slice with a label is a chain of diamonds, similar to Figure \ref{fig-chain of diamonds}, connecting a grey diamond from $\Delta_0$ to another grey diamond from $\Delta_0$.}\label{fig:diamonds}
\end{figure}
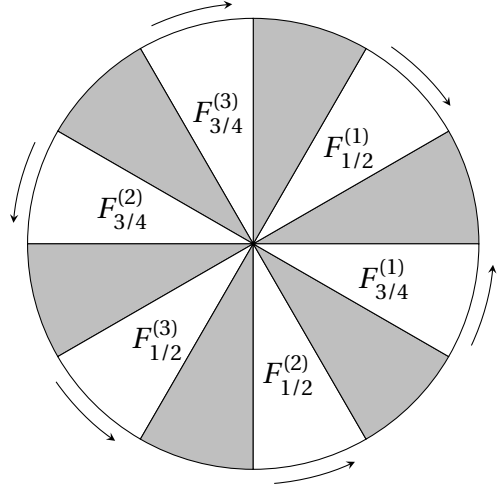

Combining the above, we conclude that there are exactly four possible such diagrams.
Indeed, up to cyclic permutation and reflections, the chains of diamonds in Figure \ref{fig:diamonds} can always be arranged so that the labels read $F_{1/2}^{(1)}F_{3/4}^{(1)}F^{(2)}_{1/2}F^{(3)}_{1/2}F^{(2)}_{3/4}F^{(3)}_{3/4}$ in clockwise order starting from the top rightmost label.
In such a sequence, we say that $F_{1/2}$ is increasing if the arrow on $F_{1/2}^{(1)}$ is clockwise and decreasing otherwise. Similarly for $F_{3/4}$ with $F_{3/4}^{(1)}$.
These choices uniquely determine the directions of the other $F^{(i)}_{1/2}$, $F^{(i)}_{3/4}$ as their directions relative to each other must remain consistent.
Consequently, there are precisely four such possibilities indexed by these choices.

As previously mentioned, the homeomorphisms $f_h$ induce homeomorphisms of the big cell and hence induce homeomorphisms of these diagrams.
One readily checks that $f_{R_1}$ takes an increasing $F_{1/2}$ (resp. decreasing $F_{1/2}$) to a decreasing $F_{1/2}$ (resp. increasing $F_{1/2}$) and an increasing $F_{3/4}$ (resp. decreasing $F_{3/4}$) to an increasing $F_{1/2}$ (resp. decreasing $F_{3/4}$).
Similarly, $f_{R_6}$ takes an increasing $F_{1/2}$ (resp. decreasing $F_{1/2}$) to an increasing $F_{1/2}$ (resp. decreasing $F_{1/2}$) and an increasing $F_{3/4}$ (resp. decreasing $F_{3/4}$) to a decreasing $F_{1/2}$ (resp. increasing $F_{3/4}$).
As we are assuming that $f_h(M_i)=M_i$, it must be the case that it preserves this diagram up to cyclic orientation and reflection. 
It is then straightforward that $h$ has this property if and only if it preserves increasing/decreasing for \emph{both} $F_{1/2}$, $F_{3/4}$. 
Analyzing $f_\mathrm{id}$, $f_{R_1}$, $f_{R_6}$, and  $f_{R_1R_6}=f_{R_1}\circ f_{R_6}$, it is apparent from the description above that this happens if and only if $f_h=f_\mathrm{id}$.
Therefore, we conclude that $h=\mathrm{id}$, as desired.
\end{proof}


In fact, examining the proof of Proposition \ref{prop:onehomeotype} one obtains the following conclusion, which proves Theorem \ref{thm:koplusoneisom}.


\begin{corollary}\label{cor:isometrieskplusone}
Suppose that $M\in\mathcal{M}_{k+1,k}$, then
$$\Isom(M)\cong\begin{cases}
\{ 1 \}, & \text{ if } \ i(M)\neq j(M),\\
\Z/2\Z, & \text{ if } \ i(M)=j(M),
\end{cases},$$
where these are the isometries with respect to the unique hyperbolic structure on $M$.
\end{corollary}


\begin{proof}
Indeed by Corollary \ref{cor:combinatorialiso}, any isometry of $M$ is induced by a combinatorial isomorphism of its triangulation.
However, when $i(M)\neq j(M)$, we saw in the proof of Proposition \ref{prop:onehomeotype} that there are no such isomorphisms as $f_h(M_i)=M_i$ for $h\in X_0$ if and only if $h=\mathrm{id}_{X_0}$.
Therefore the isometry group is trivial.

When $i(M)=j(M)$, there is an extra combinatorial automorphism, $\iota$, of $\Delta_0$ which takes $E_1$ to $E_6$ and interchanges the pairs $(F_1,F_2)$ and $(F_3,F_4)$. 
In this case, the subgroup of combinatorial automorphisms, $\widetilde{X}_0$, of $\Delta_0$ which acts on the manifolds in $\mathcal{M}$ is precisely the group $\widetilde{X}_0=\langle \iota, X_0\rangle$, which has order $8$. 
Indeed, $\iota$ extends to a combinatorial involution $f_\iota:M_i\to M_j$ via
\begin{align*}
f_{\iota}:M_i&\to M_j,\\
\Delta_0&\mapsto \iota\cdot \Delta_0,\\
M(\mathcal{K}_1)&\mapsto M(\mathcal{K}_2),\\
M(\mathcal{K}_2)&\mapsto M(\mathcal{K}_1),
\end{align*}
The arguments of Proposition \ref{prop:onehomeotype} again show that the only combinatorial isomorphism of the triangulation of $M$ contained in the subgroup $X_0$ is the identity.
By orbit stabilizer of the action on $\M$, we deduce that there is precisely one non-trivial element of $\widetilde{X}_0$, which is a combinatorial isomorphism of $M$.
This gives the desired result.
\end{proof}


\begin{remark}
As we will see in Section \ref{sec:examples}, one can computationally verify Corollary \ref{cor:isometrieskplusone} for the unique manifold in $\mathcal{M}_{2,1}$. 
Namely, one can show computationally that it has no isometries, which aligns with the statement of the corollary since it corresponds to $(i(M),j(M))=(0,1)$.
\end{remark}




\section{Dehn fillings of manifolds in \texorpdfstring{$\M_{k,k}$}{Mkk}}\label{sec:dehnfilling}

In this section, we prove Theorem \ref{thm-dehnfill} which gives a description of the Dehn fillings of manifolds in $\mathcal{M}_{k,k}$ with certain small slopes.
In order to do this, we first refine the description of the spines of manifolds $M\in\mathcal{M}_{k,k}$, using the terminology established in Section \ref{subsection:spinebackground} and in Section \ref{sec:kkeven}.

Recall that there is a unique manifold $M\in \mathcal{M}_{k,k}$ which is obtained by closing a $k$-chain $\{\Delta_1, \Delta_2, \dots, \Delta_{2k}\}$. We freely continue to use all of the notation $a_i, b_i$, $c_i$, $p_i, q_i$, $r_i$ and $A_i, B_i, C_i, P_i, Q_i$, $R_i$, as described in Section \ref{sec:kkeven} and shown in Figures \ref{fig-spine_notation}--\ref{fig-big face G}.

Then, as described in Section \ref{subsection:spinebackground}, the vertices of $P$ are $v_1, v_2,\dots, v_{2k}$ and the edges of $P$ are the pairs of spinal half edges of tetrahedra that are connected through face gluings.
For $1\leq i \leq k$, there are three edges between $v_{2i-1}$ and $v_{2i}$.
We denote these edges in the spine $P$ by $e^i_1$, $e^i_2$, and $e^i_3$.
For $1\leq i \leq k$, there is one edge between $v_{2i}$ and $v_{2i+1}$ where the subscripts are taken modulo $2k$; we denote this edge by $e^{i,i+1}$, with superscript taken modulo $k$.

\begin{figure}
\centering
\includegraphics[width=5in]{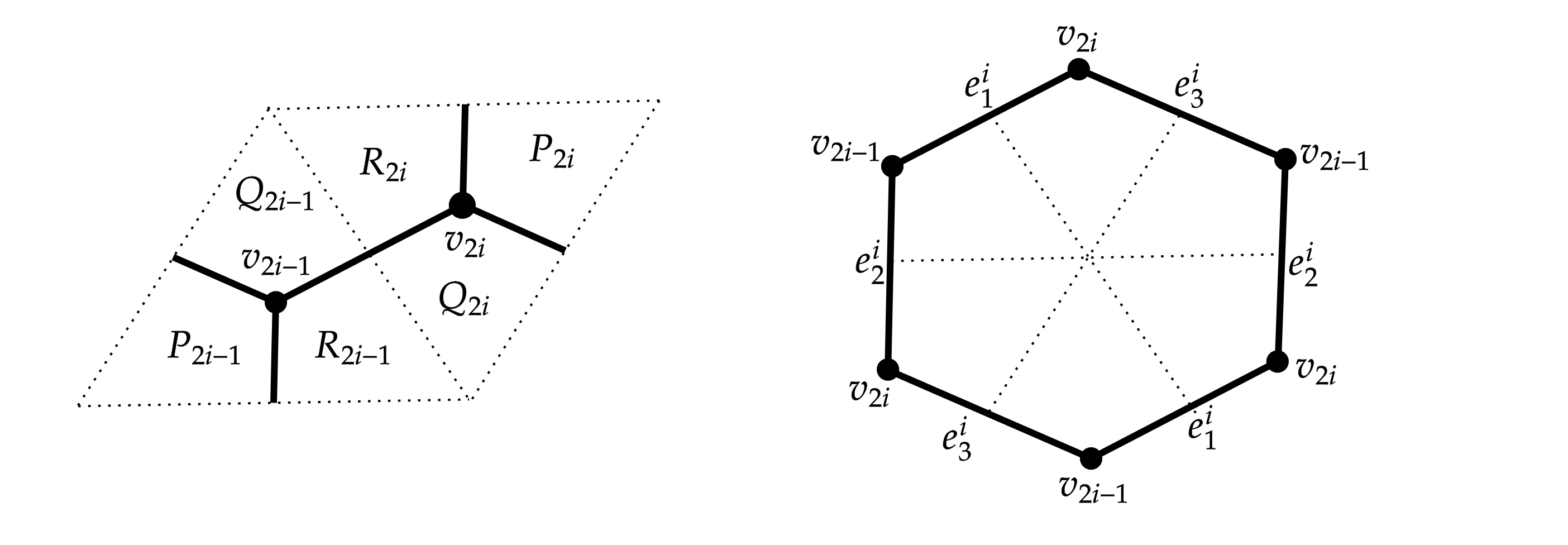}
\caption{$\Delta_{2i-1}$ and $\Delta_{2i}$ viewed from the ideal vertex they share and the hexagonal face $H_i$ of the spine}
\label{fig-hexagonal_face}
\end{figure}

The hexagonal face $H_i$ of the spine $P$ defined in Section \ref{subsection:spinebackground} is then made up of the diamonds $P_j$, $Q_j$, and $R_j$ for $j=2i-1, 2i$ (see Figure \ref{fig-hexagonal_face}).
 The remaining diamonds $A_i, B_i$, and $C_i$ for $1\leq i \leq 2k$ form the big face $G$ which was described in detail in Section \ref{sec:kkeven}.

We are now in a position to prove Theorem \ref{thm-dehnfill}, which we restate for the reader's convenience.


\dehnfill*


\begin{proof}
Let $M(\alpha)$ be the manifold obtained by Dehn filling the $i^{th}$ cusp of $M$ along a slope $\alpha$.
Let $D$ be a meridian disk of the solid torus resulting from the filling. The curve $\partial D$ cuts the hexagonal face $H_i$ of $P$ into a union of open faces of $P \cup D$.
It was shown in \cite{fmp-dehn} that the removal of one such open face from $P \cup D$ yields a spine, $P(\alpha)$, of the filled manifold $M(\alpha)$ (see \cite[\S 3]{fmp-dehn}).
We prove the theorem for $\alpha = 1/2$ and a similar argument will yield the cases for other $\alpha$.
Figure \ref{fig-dehn_fill_1/2} shows the intersection of $H_i$ with $\partial D$ when $\alpha = 1/2$.

\begin{figure}
\centering
\includegraphics[width=5in]{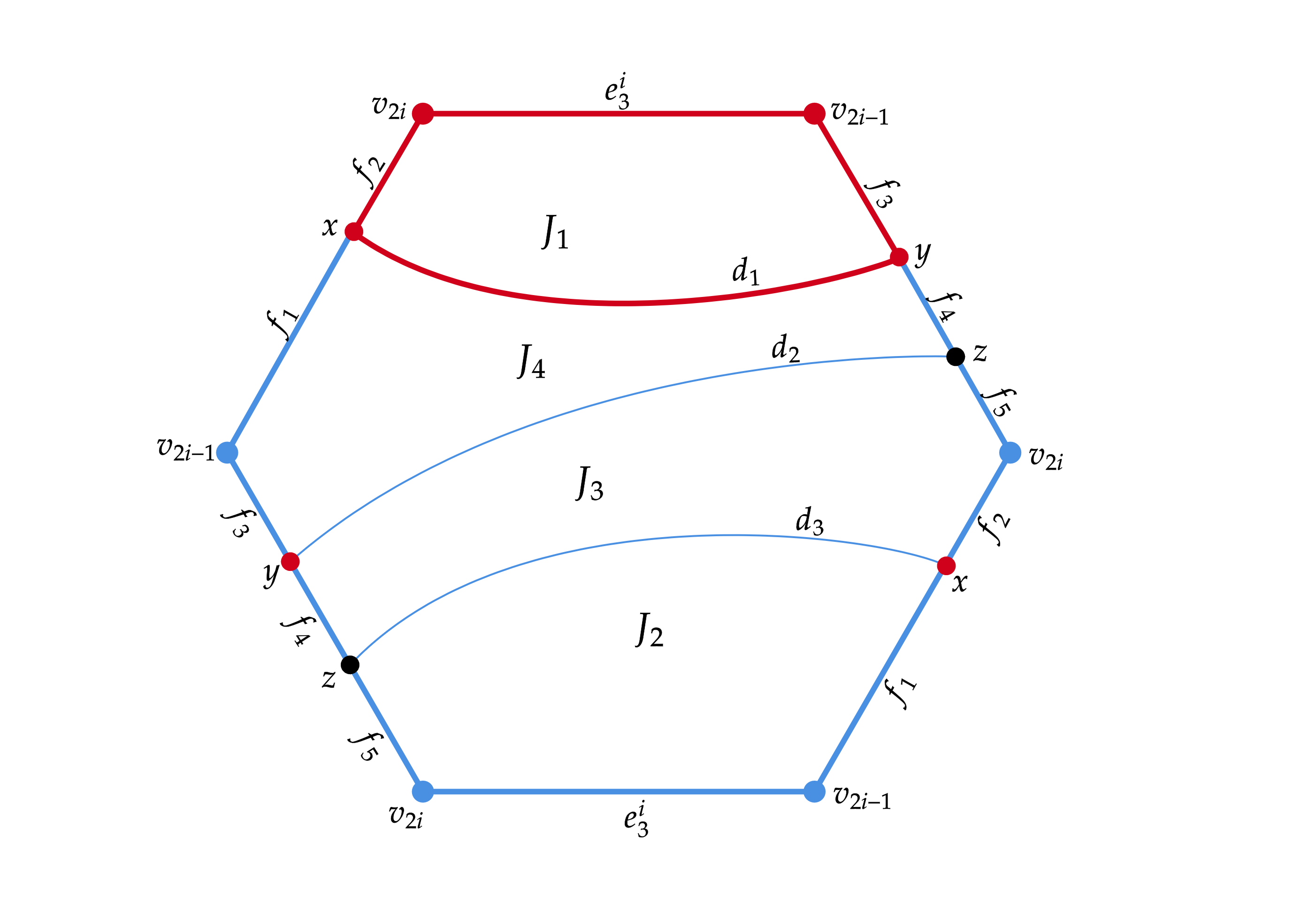}
\caption{Intersection of the slope $1/2$ with the hexagonal face $H_i$}
\label{fig-dehn_fill_1/2}
\end{figure}

We label the points at which $\partial D$ intersects the $1$-skeleton of $P$ as $x$, $y$, and $z$. These points divide $\partial D$ into three segments which we denote by $d_1$, $d_2$, and $d_3$.
Furthermore, $x$ lies on the interior of $e^i_1$ and divides $e^i_1$ into line segments $f_1$ and $f_2$. 
Similarly $y$ and $z$ divide $e^i_2$ into three segments $f_3$, $f_4$, and $f_5$ as shown in Figure \ref{fig-dehn_fill_1/2}.

If we remove the open face $J_1$, then $P(\alpha) = (P \cup D) \backslash J_1$ is a spine for $M(\alpha)=M(1/2)$. 
Below we list a preliminary cell decomposition of $P(\alpha)$ into vertices, edges, and faces denoted $V(P(\alpha))$, $E(P(\alpha))$, and $F(P(\alpha))$, respectively.
\begin{itemize}
    \item $V(P(\alpha)) =V(P) \cup \{x,y,z\}$,
    \item $E(P(\alpha)) = \bigg(E(P) \backslash \{e^i_1,e^i_2\} \bigg)\cup \left(\bigcup\limits_{j=1}^5 f_j \right)\cup\left( \bigcup\limits_{j=1}^3 d_j\right)$,
    \item $F(P(\alpha)) = \bigg(F(P)\backslash\{H_i\} \bigg)\cup\left(\bigcup\limits_{i=2}^4 J_i  \right)\cup D$.
\end{itemize}

We now update the cell decomposition of $P(\alpha)$ as described below.
Note that each edge in $E(P(\alpha))$ that was adjacent to $J_1$ is in the boundary of exactly two faces in $F(P(\alpha))$: one of $G$, $D$ and one of $J_2$, $J_3$, $J_4$. 
We first delete these four edges from the preliminary cell decomposition, merging $J_2$, $J_3$, $J_4$, $D$, and $G$ into a single face $G'$.

\begin{figure}
\centering
\includegraphics[width=5in]{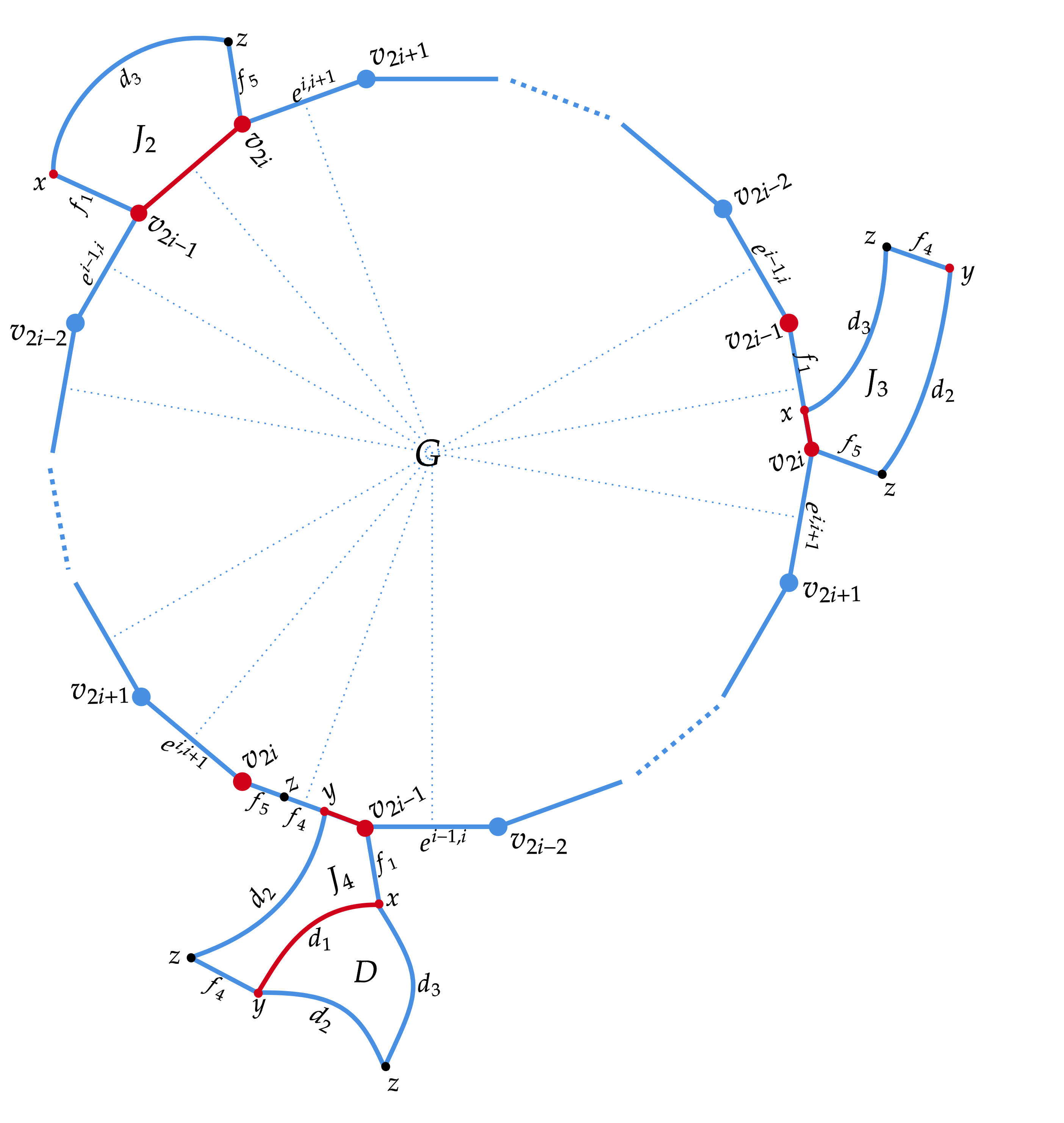}
\caption{Merging $J_2$, $J_3$, $J_4$, $D$, and $G$ into a single $2$-cell $G'$. The edges and vertices in red are then deleted from the cell decomposition.}
\end{figure}
 
Next, we delete the vertices $v_{2i-1}$, $v_{2i}$, $x$, and $y$ and merge any two edges in $E(P(\alpha))$ into a single edge if they share a deleted vertex.
In particular, this process will merge each subset $\{e^{i,i+1}, f_5\}$, $\{f_4,d_2\}$, and $\{e^{i-1,i},f_1,d_3\}$ of $E(P(\alpha))$ into three single edge classes that have multiplicity $3$ in the boundary of $G'$.
We denote these newly formed edge classes by $e^0_1$, $e^0_2$, and $e^0_3$.  

\begin{figure}[t]
\centering
\includegraphics[width=5in]{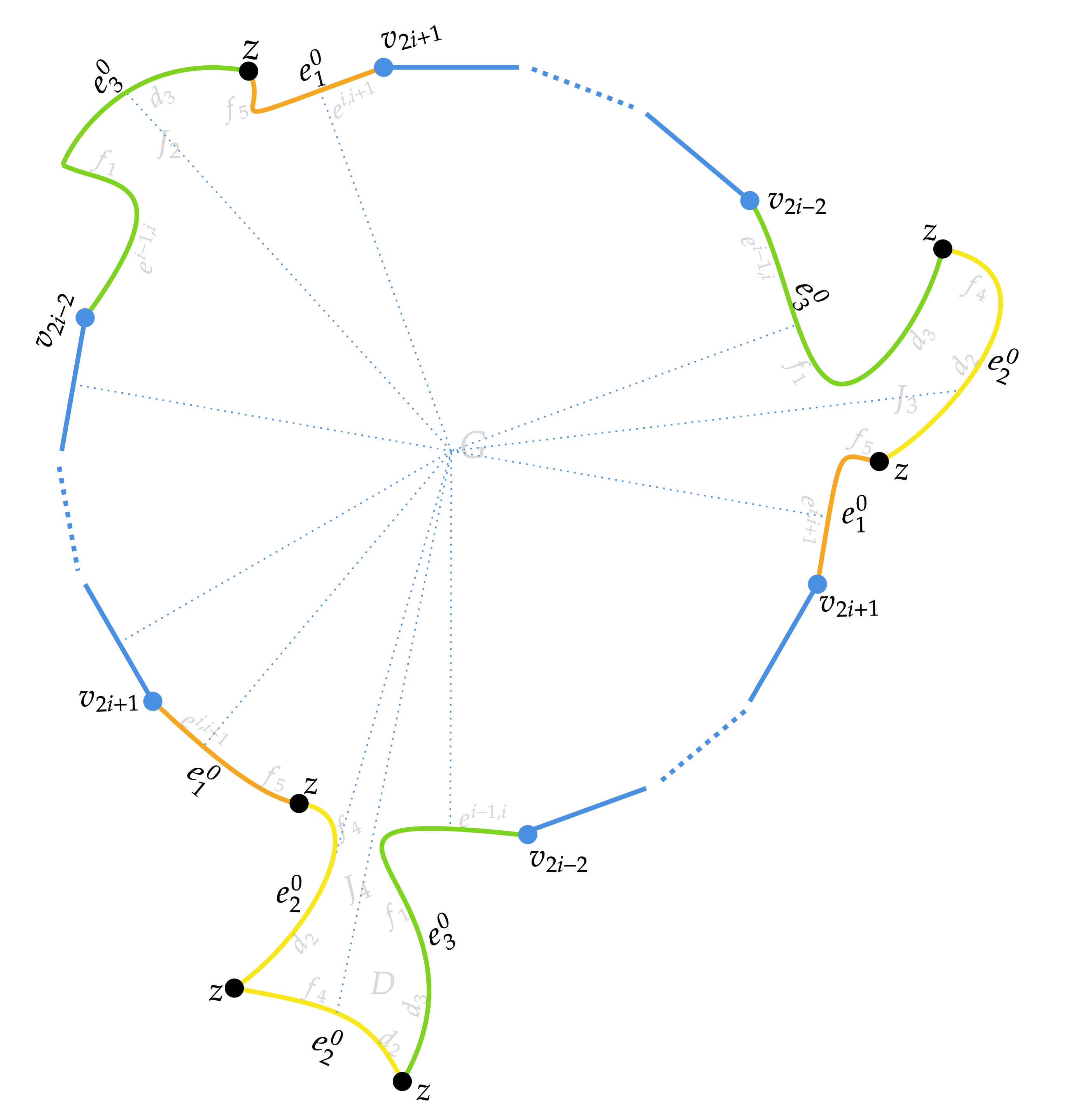}
\caption{$G'$ with its new cell decomposition. The new edge classes are $e^0_1$ in orange, $e^0_2$ in yellow, and $e^0_3$ in green. Each edge has valence $3$ in $\partial G'$.}
\end{figure}

Note that this process does not alter the cell decomposition of hexagonal faces $H_j$ for which $j\neq i$, as no such face contains any of the deleted vertices or edges in its boundary.
The updated cell decomposition of $P(\alpha)$ after this process is as follows:
\begin{itemize}
    \item $\widetilde{V}(P(\alpha)) = \{v_j \mid 1\leq j \leq 2k , j\neq 2i-1,2i\} \cup \{z\}$,
    \item $\widetilde{E}(P(\alpha)) = \{e^j_1,e^j_2,e^j_3 \mid 1\leq j \leq k , j\neq i\} \cup \left\{e^{j,j+1} \mid j\neq i ,i-1 \right\} \cup \left\{e^0_1,e^0_2,e^0_3\right\}$,
    \item $\widetilde{F}(P(\alpha)) = \{H_j \mid 1\leq j \leq k , j\neq i\} \cup G'$.
\end{itemize}
It is then straightforward to see $P(\alpha)$ with this cell decomposition coincides with the spine of the unique member in $\mathcal{M}_{k,k-1}$ corresponding to the trivial partition of its set of boundary components, as described at the end of the proof of Proposition \ref{prop:onehomeotype}. 

In particular, the ideal triangulation $\mathcal{T}$ that is dual to the new cell decomposition of $P(\alpha)$ consists of a compact, fully truncated tetrahedron $\Delta_0$ (that is dual to the vertex $z$ of $P(\alpha)$) and the $k-1$ chain $\{\Delta_{2i+1}, \Delta_{2i+2}, \dots , ,\Delta_{2k-1},\Delta_{2k},\Delta_{1},\Delta_{2},\dots \Delta_{2i-3},\Delta_{2i-2}\}$. 
The boundary faces of this $k-1$ chain are the finite faces of $\Delta_{2i+1}$ and $\Delta_{2i-1}$, each of which is identified with a compact face of $\Delta_0$ in $\mathcal{T}$. 
These faces of $\Delta_0$ are dual to $e^0_1$ and $e^0_3$ in $\widetilde{E}(P(\alpha))$.
The remaining two faces of $\Delta_0$ are identified with each other and become dual to $e^0_2$ in $\widetilde{E}(P(\alpha))$.
\end{proof}




\section{Examples}\label{sec:examples}


\subsection{Topological invariants of the manifolds in \texorpdfstring{$\mathcal{M}_{k,k}$}{Mkk}}
In this subsection, we explore topological invariants of manifolds in $\mathcal{M}_{k,k}$ with the ultimate goal of proving the following theorem, which we restate for the reader's convenience.


\kevenquasiarithmetic*


Before giving the proof, we first remind the reader of some well-known facts about Coxeter polyhedra as well as the relevant definitions of the topological invariants from this theorem.
We follow the exposition of \cite[\S 10.4]{MRBook}, which details the work of Vinberg \cite{Vinberg}.

To this end, let $\mathcal{P}$ be a polyhedron in $\Hy^3$ with totally geodesic faces whose dihedral angles are submultiples of $\pi$.
Let $\Gamma_\mathcal{P}<\Isom(\Hy^3)$ be the group generated by reflections in the faces of $\mathcal{P}$. Then $\Gamma_\mathcal{P}$ is a discrete subgroup whose (index $2$) orientation preserving subgroup, $\Gamma^+_\mathcal{P}$, is a discrete subgroup of $\Isom^+(\Hy^3)$.
When $\mathcal{P}$ is compact (resp. finite volume) it follows that $\Gamma^+_\mathcal{P}$ is a cocompact (resp. cofinite) lattice.
In the case that $\Gamma^+_\mathcal{P}$ is a lattice, we may describe its topological invariants using a procedure discovered by Vinberg \cite{Vinberg} and detailed below.

Enumerate the faces of $\mathcal{P}$ by $F_1,\dots, F_n$ and let $\{\vec{e}_1,\dots,\vec{e}_n\}$ be the corresponding outward pointing unit normals to the faces $F_i$ in the hyperboloid model of $\Hy^3$.
Then the \emph{Gram matrix} is the matrix $\mathcal{G}_\mathcal{P}=(a_{ij})$ whose coefficients are $a_{ij}=B(\vec{e_i},\vec{e_j})$, where $B(-,-)$ is the standard symmetric bilinear form giving rise to $\Hy^3\subset \R^4$.
More, explicitly
$$a_{ij}=\begin{cases}
2,&i=j,\\
-2\cos(\theta_{ij}),& F_i\cap F_j\neq\emptyset,~i\neq j\\
-2\cosh(\ell_{ij}),& F_i\cap F_j=\emptyset,
\end{cases},$$
where $\theta_{ij}$ is the dihedral angle between $F_i$, $F_j$ when they intersect and $\ell_{ij}$ is the hyperbolic length between $F_i$, $F_j$ when they do not.
Given an ordered subset $I=\{i_1,\dots,i_r\}\subseteq\{1,\dots, n\}$ we define a \emph{cyclic product} by the formula $b_I=a_{i_1i_2}a_{i_2i_3}\cdots a_{i_{r-1}i_r}a_{i_ri_1}$ and define the set of all cyclic products by
$$\Pi_{\mathcal{G}_\mathcal{P}}=\{b_I\mid I\subseteq\{1,\dots,n\}\}.$$
It is a result of Vinberg \cite{Vinberg} that the field generated by $\Pi_{\mathcal{G}_\mathcal{P}}$ over $\Q$, denoted $\Q(\Pi_{\mathcal{G}_\mathcal{P}})$, is precisely the adjoint trace field of $\Gamma^+_\mathcal{P}$.

\begin{remark}
There is a minor discrepancy in the description of the Gram matrix between \cite{MRBook} and \cite{Vinberg}. 
Namely, \cite{MRBook} defines the entries of the Gram matrix to be twice those of \cite{Vinberg}. 
Though this discrepancy seems innocuous, the factor of $2$ becomes relevant when performing the calculation for quasi-arithmeticity described below.
We remind the reader that we are following the notation of \cite{MRBook} and therefore they must take this into account when matching these statements to \cite{Vinberg}.
\end{remark}

In this paper, we are interested in the invariant trace field of $\Gamma^+_\mathcal{P}$ in the sense of \cite{MRBook}, which is a quadratic extension of the adjoint trace field.
We now describe how to obtain the invariant trace field from the adjoint trace field in this setting. 
For any ordered subset $I=\{i_1,\dots,i_r\}$ as above, define
$$\vec{v}_I=a_{i_1i_2}a_{i_2i_3}\cdots a_{i_{r-1}i_r}\vec{e}_{i_r},$$
with $\vec{e}_{i_r}$ as before.
The $\Q(\Pi_{\mathcal{G}_\mathcal{P}})$-span of the vectors $\vec{v}_I$ as $I$ varies over all ordered subsets of $\{1,\dots, n\}$ is a 4-dimensional $\Q(\Pi_{\mathcal{G}_\mathcal{P}})$-quadratic space $(W,q)$ such that $W\otimes \R\cong \R^4$ has signature $(3,1)$. 
Letting $d\in \Q(\Pi_{\mathcal{G}_\mathcal{P}})^*$ be any lift of the discriminant $\mathrm{disc}(q)\in \Q(\Pi_{\mathcal{G}_\mathcal{P}})^*/\Q(\Pi_{\mathcal{G}_\mathcal{P}})^{*2}$, it follows from \cite[Thm 10.4.1]{MRBook} that the invariant trace field of $\Gamma^+_\mathcal{P}$ is given by $\Q(\Pi_{\mathcal{G}_\mathcal{P}})(\sqrt{d})$. 
It is clear from the definition that the latter is independent of the choice of lift.

Let $\Gamma<\Isom^+(\Hy^3)$ be a lattice which is not cocompact. 
Then we say that $\Gamma$ is \emph{quasi-arithmetic} if the invariant trace field is imaginary quadratic\footnote{There is a more general definition of a quasi-arithmetic lattice but this will suffice for our purposes, since all manifolds in this paper are finite-volume, noncompact.}.
Additionally, we say that a lattice $\Gamma<\Isom^+(\Hy^3)\cong\mathrm{PSL}_2(\C)$ has \emph{integral traces} if the set of (lifts of) traces of $\Gamma$, $\{\mathrm{Tr}(\widetilde\gamma^2)\mid \gamma\in\Gamma\}$, is contained in the ring of integers of its invariant trace field.
Here the notation $\widetilde\gamma$ denotes any choice of lift of $\gamma$ to $\mathrm{SL}_2(\C)$.
We also note that this property is independent of choice of lift and, necessarily, $\{\mathrm{Tr}(\widetilde\gamma^2)\mid \gamma\in\Gamma\}$, is always contained in the invariant trace field (though not necessarily its ring of integers).
If $\Gamma$ is both quasi-arithmetic and has integral traces, then it is arithmetic and commensurable with a Bianchi group $\mathrm{PSL}_2(\mathcal{O}_d)$, where $\mathcal{O}_d$ is the ring of integers of an imaginary quadratic extension of $\Q$.

In the setting that $\Gamma=\Gamma^+_{\mathcal{G}_\mathcal{P}}$, Vinberg shows that $\Gamma$ has integral traces if and only if $\Pi_{\mathcal{G}_\mathcal{P}}\subset\mathcal{O}_{\Q(\Pi_{\mathcal{G}_\mathcal{P}})}$, the ring of integers of $\Q(\Pi_{\mathcal{G}_\mathcal{P}})$, which is also in turn equivalent to the entries of $\mathcal{G}_\mathcal{P}$ being algebraic integers.
With this in mind, we are now in a position to prove Theorem \ref{thm:kevenquasiarithmetic}.


\begin{figure}[t]
\begin{overpic}[width=1.75in]{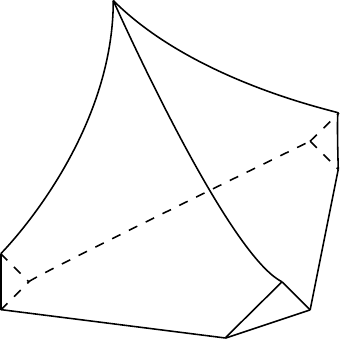}
\put(17, 70){\tiny$ \pi/3$}
\put(60, 83){\tiny$ \pi/3$}
\put(49, 70){\tiny$ \pi/3$}
\put(40, 38.5){\tiny$ \theta_k$}
\put(28, -1.5){\tiny$ \theta_k$}
\put(98, 30){\tiny$ \theta_k$}
\end{overpic}
\caption{The geometric realization of a truncated tetrahedron with a single ideal vertex. All angles which are not labeled are right angles and all faces are subsets of totally geodesic copies of $\Hy^2$ in $\Hy^3$.}\label{fig:truncatedtetrawithangles}
\end{figure}

\begin{proof}[Proof of Theorem \ref{thm:kevenquasiarithmetic}]
Fix some $k$ even and let $\theta_k=\pi/3k$.
Let $\mathcal{P}_k$ be the truncated tetrahedron from Figure \ref{fig:truncatedtetrawithangles} with dihedral angles labeled therein.
Define the following three quantities
$$
z_k = -2\cos(\theta_k),\quad\quad
\ell_1 = \cosh^{-1}\left( \frac{3\cos(\theta_k)}{\sqrt{1 + 2\cos(2\theta_k)}}\right),\quad\quad
\ell_2 = \cosh^{-1}\left( \frac{1/2 + \cos^2(\theta_k)}{1/2 + \cos(2\theta_k)}\right),
$$
which are proportional to entries of the Gram matrix for which two faces intersect with dihedral angle $\theta_k$, an original face and a triangular face do not intersect, and two triangular faces do not intersect, respectively. 
In particular, the length $\ell_2$ is the length of the edge labeled with $\theta_k$ in Figure \ref{fig:truncatedtetrawithangles}.
The value for $z_k$ follows by definition and we include the derivation of $\ell_1$, $\ell_2$ in what follows below.
The reader willing to take these quantities on faith can skip to the description of the Gram matrix in this proof.

\begin{figure}[t]
\includegraphics[width=4in]{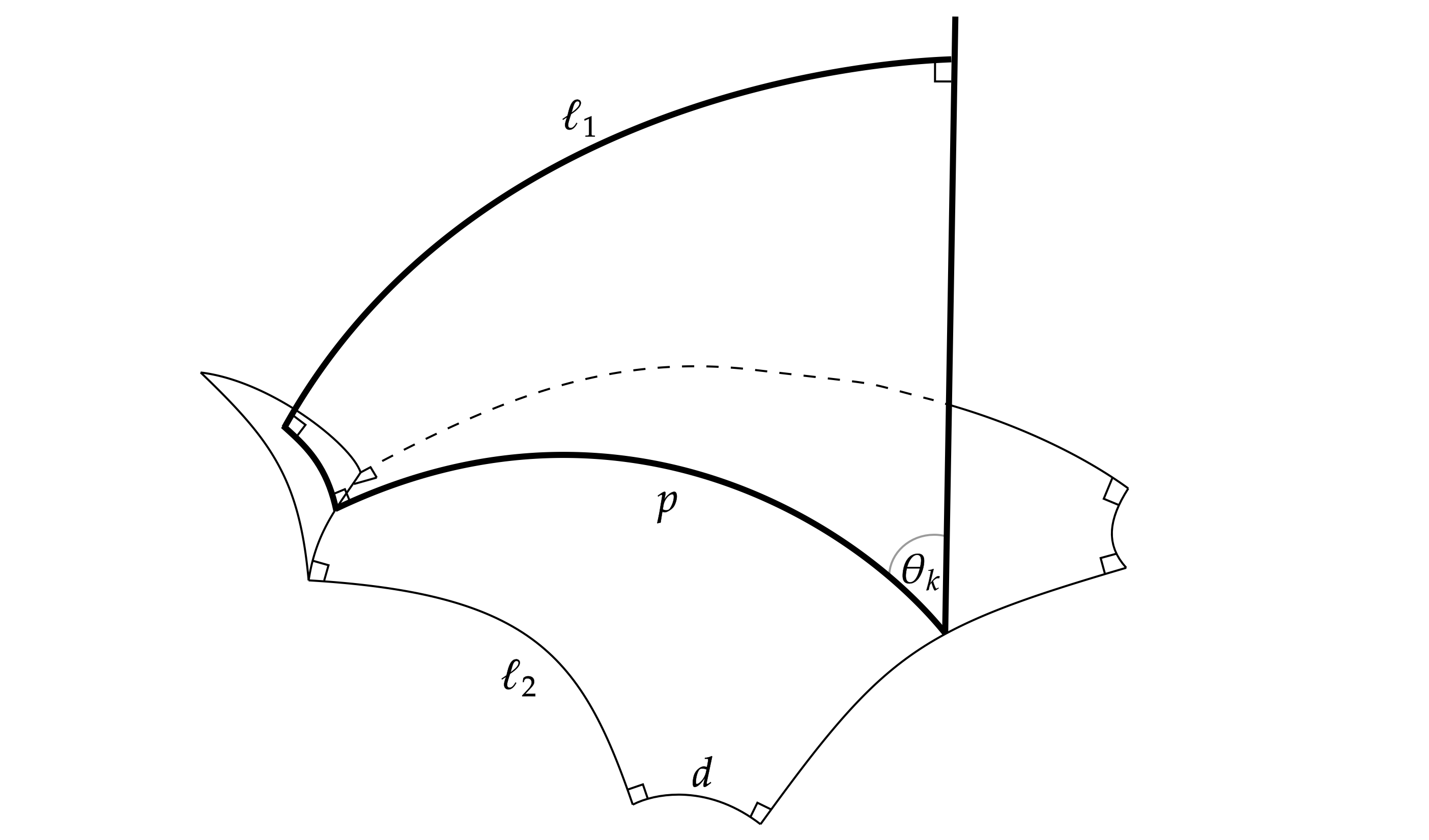}
\caption{The polygons one uses to compute quantities $\ell_1$ and $\ell_2$.}\label{fig:trigonometry}
\end{figure}

Any of the three triangles arising as triangular faces have angles $\pi/3$, $\theta_k$, and $\theta_k$.
The triangle intersects the finite face of $\mathcal{P}_k$ in the side opposite to the angle $\pi/3$.
Let $d$ be the length of this side.  
By the second law of cosines \cite[Thm 3.5.4]{ratcliffe1994foundations}, we have
\[\cosh(d)=\frac{1/2+ \cos^2(\theta_k)}{\sin^2(\theta_k)}.\]
The finite face of $\mathcal{P}_k$ is a right-angled hexagon and has three sides of length $d$, which alternate every other side (see Figure \ref{fig:trigonometry}). By the hexagon rule \cite[Thm 3.5.14]{ratcliffe1994foundations}, the length of any of the remaining sides can be computed by 
\[\cosh(\ell_2)= \frac{\cosh^2(d)+\cosh(d)}{\sinh^2(d)} = 1 + \frac{1}{\cosh(d)-1}.\]
By replacing $\cosh(d)$ with $\frac{1/2+ \cos^2(\theta_k)}{\sin^2(\theta_k)}$, we obtain that
\[\cosh(\ell_2)=\frac{1/2+\cos^2(\theta_k)}{1/2+\cos(2\theta_k)}.\]

Now consider a side of length $d$ in the finite face and the side of length $\ell_2$ opposite to it. Let $p$ be the length of the common perpendicular to these two sides. This common perpendicular divides the finite face into two isometric right-angled pentagons (see Figure \ref{fig:trigonometry}). By \cite[Thm 3.5.11]{ratcliffe1994foundations} we see that
\[\cosh(p)= \tanh(d/2)\tanh(\ell_2/2).\]
Replacing $d$ and $\ell_2$ using the first two equations, we obtain
\[\cosh(p)= \frac{3\cos(\theta_k)}{\sin(\theta_k)\sqrt{2\cos^2(\theta_k)+1}}.\]

Finally, we consider the pair of faces of $\mathcal{P}_k$ given by an original face of $\mathcal{P}_k$ that is adjacent to its ideal vertex and the triangular face that is opposite to it. Let $\ell_1$ be the length of the common perpendicular to these two faces. Then there is a Lambert quadrilateral (a quadrilateral with three right angles) with opposite side lengths of $p$ and $\ell_2$. The fourth angle of this quadrilateral is $\theta_k$ and is adjacent to the side with length $p$. By applying a law of cosines \cite[Thm 3.5.7]{ratcliffe1994foundations} we then obtain 
\[\cosh(\ell_1) = \sin(\theta_k) \cosh(p).\]
This confirms the formulas above.

Define $\widetilde{\ell}_i=-2\cosh(\ell_i)$ for $i=1,2$, then the Gram matrix of the corresponding truncated tetrahedron is given by
$$G_k =
\begin{pmatrix}
2 & -1 & -1 & z_k & 0 & 0 & \widetilde{\ell}_1 \\
-1 & 2 & -1 & z_k & \widetilde{\ell}_1 & 0 & 0 \\
-1 & -1 & 2 & z_k & 0 & \widetilde{\ell}_1 & 0 \\
z_k & z_k & z_k & 2 & 0 & 0 & 0 \\
0 & \widetilde{\ell}_1 & 0 & 0 & 2 & \widetilde{\ell}_2 & \widetilde{\ell}_2 \\
0 & 0 & \widetilde{\ell}_1 & 0 & \widetilde{\ell}_2 & 2 & \widetilde{\ell}_2 \\
\widetilde{\ell}_1 & 0 & 0 & 0 & \widetilde{\ell}_2 & \widetilde{\ell}_2 & 2
\end{pmatrix}.$$
It is straightforward to check, using Mathematica for instance, that the set of distinct cyclic products of $G_k$ is given by
$$\Pi_k=\left\{2, 1, z_k^2, \widetilde{\ell}_1~\!\!\!\!\!^2, \widetilde{\ell}_2~\!\!\!\!\!^2, -1, -z_k~\!\!\!\!\!^2, \widetilde{\ell}_2~\!\!\!\!\!^3, -\widetilde{\ell}_1~\!\!\!\!\!^2 \widetilde{\ell}_2, 
 \widetilde{\ell}_1~\!\!\!\!\!^2 \widetilde{\ell}_2, -\widetilde{\ell}_1~\!\!\!\!\!^2 \widetilde{\ell}_2~\!\!\!\!\!^2\right\}.$$
Since $\widetilde{\ell}_1~\!\!\!\!\!^2$, $\widetilde{\ell}_2$ are rational functions in $\cos^2(\theta_k)$, one readily sees that
$$N_k=\Q(\Pi_k)=\Q(z_k^2,\widetilde{\ell}~\!^2_1,\widetilde{\ell}_2)=\Q\left(\cos^2\left(\frac{\pi}{3k}\right)\right)=\Q\left(\cos\left(\frac{2\pi}{3k}\right)\right),$$
where the last line is the standard trig identity $\cos^2(x)=\frac{1}{2}\left(1+\cos(2x)\right)$.
We remark that the adjoint trace field, $N_k$, is always totally real and is equal to $\Q$ if and only if $k=2$.

Following the procedure outlined at the beginning of this section, we may define $\vec{v}_1=2\vec{e}_1$, $\vec{v}_2=-\vec{e}_2$, $\vec{v}_3=-\vec{e}_3$, $\vec{v}_4=z_k\vec{e}_4$, and restrict the Gram matrix to $W=\{\vec{v}_1,\dots,\vec{v}_4\}$ which yields a $4$-dimensional $N_k$-quadratic space with associated bilinear form given by
$$G'_k= \begin{pmatrix}
8 & 2 & 2 & 2z_k^2 \\
2 & 2 & -1 & -z_k^2 \\
2 & -1 & 2 & -z_k^2 \\
2z_k^2 & -z_k^2 & -z_k^2 & 2z_k^2
\end{pmatrix}.$$
Again, using Mathematica for instance, one checks that $\det(G'_k)=-108 z_k^4$ and therefore the discriminant of the associated quadratic form is given by the square class $\mathrm{\mathrm{disc}(G'_k)}=[-3]\in N_k^*/N_k^{*2}$.
Consequently, the result of Vinberg shows that the invariant trace field of $\Gamma^+_{\mathcal{P}_k}$ is given by $L_k=N_k\left(\sqrt{\mathrm{disc}(G'_k)}\right)=N_k(\sqrt{-3})$.
As $N_k$ is totally real, $L_k$ is imaginary quadratic if and only if $N_k$ is equal to $\Q$ which occurs precisely when $k=2$.
Therefore $\Gamma^+_{\mathcal{P}_k}$ is quasi-arithmetic if and only if $k=2$.

We now turn to the statement on integrality. 
As $z^2_k\in\mathcal{O}_{N_k}$ for all $k$, we note that if $\widetilde{\ell}~\!\!^2_1,\widetilde{\ell}_2\in\mathcal{O}_{N_k}$ then it is clear that $\Pi_k\subset \mathcal{O}_{N_k}$ and therefore that $\Gamma^+_{\mathcal{P}_k}$ has integral traces.
We will show that $\widetilde{\ell}~\!^2_1,\widetilde{\ell}_2\in\mathcal{O}_{N_k}$ if and only if $k$ is not a power of $2$ and therefore deduce, using the previous remark, that $\Gamma^+_{\mathcal{P}_k}$ has integral traces if and only if $k$ is not a power of $2$.

To this end, note that if $\alpha_k=2\cos(\theta_k)$ then 
$$\widetilde{\ell}~\!^2_1=\frac{9\alpha_k^2}{\alpha_k^2-1}=9\left(1+\frac{1}{\alpha_k^2-1}\right),\quad \widetilde{\ell}_2=-\frac{\alpha_k^2+2}{\alpha_k^2-1}=-\left(1+\frac{3}{\alpha_k^2-1}\right).$$
Consequently $\widetilde{\ell}~\!^2_1\in\mathcal{O}_{N_k}$ if and only if $9/(\alpha_k^2-1)\in \mathcal{O}_{N_k}$ and similarly $\widetilde{\ell}_2\in\mathcal{O}_{N_k}$ if and only if $3/(\alpha_k^2-1)\in \mathcal{O}_{N_k}$.
We claim that
$$\textrm{Nrm}_{N_k/\Q}(\alpha_k^2-1)=\begin{cases}
\pm 2,&k=2^e\\
\pm 1,&k\neq 2^e
\end{cases},$$
which shows that $\widetilde{\ell}~\!^2_1,\widetilde{\ell}_2\in\mathcal{O}_{N_k}$ if and only if $k$ is not a power of $2$.
Indeed, if $k$ is a power of $2$ then a norm calculation shows that $3/(\alpha_k^2-1)\notin \mathcal{O}_{N_k}$ and if $k$ is not a power of $2$ then $\alpha_k^2-1$ is a unit in $\mathcal{O}_{N_k}$ and hence $3/(\alpha_k^2-1)\in \mathcal{O}_{N_k}$.

To compute this norm, one checks that
$$\alpha_k^2-1=\eta^2_{6k}+\frac{1}{\eta^2_{6k}}+1=\frac{1}{\eta^2_{6k}}\left(1+\eta^2_{6k}+\eta_{6k}^4\right)=\frac{1}{\eta^2_{6k}}\Phi_{3}\left(\eta^2_{6k}\right),$$ where $\eta_{6k}$ is a primitive $6k^{th}$ root of unity and $\Phi_3(x)=1+x+x^2$ is the $3^{rd}$ cyclotomic polynomial.
Since the norm is multiplicative and $\eta_{6k}$ is a root of unity,
$$\textrm{Nrm}_{N_k/\Q}(\alpha_k^2-1)=\textrm{Nrm}_{N_k/\Q}\left(\frac{1}{\eta^2_{6k}}\right)\textrm{Nrm}_{N_k/\Q}\left(\Phi_{3}(\eta^2_{6k})\right)=\textrm{Nrm}_{N_k/\Q}\left(\Phi_{3}(\eta^2_{6k})\right),$$
therefore it suffices to compute the final quantity.
For this, we will instead compute the norm from the full cyclotomic field, that is, we will show that
$$\textrm{Nrm}_{\Q(\eta_{6k})/\Q}\left(\Phi_{3}(\eta^2_{6k})\right)=\begin{cases}
 2^4,&k=2^e\\
 1,&k\neq 2^e
\end{cases},$$
and then using the standard fact that
$$\textrm{Nrm}_{\Q(\eta_{6k})/\Q}(x)=\left(\textrm{Nrm}_{N_k/\Q}(x)\right)^{[\Q(\eta_{6k}):N_k]}=\left(\textrm{Nrm}_{N_k/\Q}(x)\right)^4,$$
for any $x\in N_k$, we will conclude our claim.
For the claimed computation, recall that
$$\textrm{Nrm}_{\Q(\eta_{6k})/\Q}\left(\Phi_{3}(\eta^2_{6k})\right)=\prod_{\substack{i=1 \\ \gcd(i,6k)=1}}^{6k}\left(1+\eta_{6k}^{2i}+\eta_{6k}^{4i}\right)=\prod_{\substack{i=1 \\ \gcd(i,6k)=1}}^{6k}\Phi_{3}(\eta^{2i}_{6k}).$$
As $\gcd(i,6k)=1$, each $\eta_{6k}^i$ is also a primitive $6k^{th}$ root of unity and consequently, as we vary over this set, the $\eta^{i}_{6k}$ comprise all roots of the cyclotomic polynomial $\Phi_{6k}(x)$.
In particular, this means that the latter quantity is the resultant of $\Phi_{6k}(x)$ and $\Phi_3(x^2)$.
Note that
$$\Phi_3(x^2)=1+x^2+x^4=(1+x+x^2)(1-x+x^2)=\Phi_3(x)\Phi_6(x),$$
and therefore multiplicativity of the resultant shows that
$$\textrm{Nrm}_{\Q(\eta_{6k})/\Q}\left(\Phi_{3}(\eta^2_{6k})\right)=\mathrm{Res}(\Phi_{6k}(x),\Phi_3(x^2))=\mathrm{Res}(\Phi_{6k}(x),\Phi_3(x))\cdot\mathrm{Res}(\Phi_{6k}(x),\Phi_6(x)).$$
Resultants of cyclotomic polynomials were computed by Apostol \cite{Apostol}, who showed that if $m>n$ then $\mathrm{Res}(\Phi_m(x),\Phi_n(x))=1$ if and only if $m/n$ is not a prime power.
Recall that $k$ is even and therefore $6k/3=2k$ and $6k/6=k$ are prime powers if and only if $k$ is $2^e$ for some $e\in\N$. Moreover, the calculation of \cite{Apostol} shows that
$$\mathrm{Res}(\Phi_{6k}(x),\Phi_3(x))\cdot\mathrm{Res}(\Phi_{6k}(x),\Phi_6(x))=\begin{cases}
 2^{\varphi(3)}2^{\varphi(6)},&k=2^e\\
 1,&k\neq 2^e
\end{cases}=\begin{cases}
 2^4,&k=2^e\\
 1,&k\neq 2^e
\end{cases},$$
where $\varphi(-)$ is Euler's phi function.
This calculation completes the proof of integrality.

Finally, we claim that $\Gamma^+_{\mathcal{P}_k}$ is commensurable with the lattice $\pi_1(DM)<\Isom^+(\Hy^3)$.
This will complete the proof as the invariant trace field, quasi-arithmeticity, and non-integrality are commensurability invariants and the invariant quaternion algebras of finite-volume, non-compact hyperbolic $3$-manifolds are always matrix algebras over the invariant trace field (see \cite[Thm 3.3.8]{MRBook}).

To this end, note that $DM$ can be decomposed into a union of geometric tetrahedra isometric to $\mathcal{P}_k$.
Therefore, on lifting to the universal cover, we obtain a tiling $\mathcal{T}(\mathcal{P}_k)$ of $\Hy^3$ by truncated geometric tetrahedra and so $\pi_1(DM)$ is a finite index subgroup of the group of isometries, $\Isom(\mathcal{T}(\mathcal{P}_k))$, of this tiling.
The latter is also a discrete subgroup of $\Isom(\Hy^3)$.
As the Coxeter group $\Gamma^+_{\mathcal{P}_k}$ similarly acts by isometries on $\mathcal{T}(\mathcal{P}_k)$, it also is contained in $\Isom(\mathcal{T}(\mathcal{P}_k))$ as a subgroup of finite index and therefore $\pi_1(DM)$, $\Gamma^+_{\mathcal{P}_k}$ are commensurable as required. This completes the proof.
\end{proof}


\subsection{Volumes of manifolds in \texorpdfstring{$\M_{k,k}$}{Mkk}}\label{subsection:volume}


In this subsection, we compute the volume of the unique hyperbolic $3$-manifold in $M_{k,k}$ when $k$ is even.
We retain all notation from the previous subsection, especially the definition of the angle $\theta_k$.

To this end, Ushijima \cite[Thm 1.1]{Ushijima} gave a formula for the computation of the volume of hyperbolic geometric partially truncated tetrahedra using only their dihedral angles, which we now recall. 
If $\theta_{ij}$ for $1\le i< j\le 4$ denotes the dihedral angle between the original faces $F_1$, $F_2$, $F_3$, $F_4$, then let
$$ a = e^{i\theta_{12}},\quad b = e^{i\theta_{13}},\quad c = e^{i\theta_{14}},\quad d = e^{i\theta_{23}},\quad e = e^{i\theta_{24}},\quad f = e^{i\theta_{34}}.$$
Using these, one define the quantities
\begin{align}
Z_1 &=\frac{-2\sin(\theta_{12})\sin(\theta_{34}) +\sin(\theta_{13})\sin(\theta_{24}) +\sin(\theta_{14})\sin(\theta_{23}) -\sqrt{\det G}}{ad + be + cf + abf + ace + bcd + def + abcdef},\label{eqn:Z1}\\
 Z_2 &=\frac{-2\sin(\theta_{12})\sin(\theta_{34}) +\sin(\theta_{13})\sin(\theta_{24}) +\sin(\theta_{14})\sin(\theta_{23}) +\sqrt{\det G}}{ad + be + cf + abf + ace + bcd + def + abcdef}\label{eqn:Z2},
\end{align}
where $G$ is the $4\times 4$ Gram matrix associated to the faces $F_1,\dots, F_4$.
Then the theorem of Ushijima states that the volume of the associated truncated tetrahedron is given by
\begin{equation}\label{eqn:volumetetra}
V=V(\{\theta_{ij}\})=\frac{1}{2}\Im\bigl(U(Z_1) - U(Z_2)\bigr),
\end{equation}
where
\begin{align} U(Z_k) = \tfrac{1}{2}\Bigl[&\mathrm{Li}_2\bigl(Z_k\bigr)+\mathrm{Li}_2\bigl(abdeZ_k\bigr)+\mathrm{Li}_2\bigl(acdfZ_k\bigr)+\mathrm{Li}_2\bigl(bcefZ_k\bigr)\label{eqn:Uz}\\&-\mathrm{Li}_2\bigl(abcZ_k\bigr)-\mathrm{Li}_2\bigl(-aefZ_k\bigr)-\mathrm{Li}_2\bigl(-bdfZ_k\bigr)-\mathrm{Li}_2\bigl(-cdeZ_k\bigr)\Bigr] \nonumber,
\end{align}
and $\mathrm{Li}_2(-)$ is Spence's dilogarithm.

When $M_k$ is the unique manifold in $\M_{k,k}$, we have $\theta_{12}=\theta_{13}=\theta_{23}=\pi/3$ and $\theta_{14}=\theta_{24}=\theta_{34}=\theta_k$.
Using the fact that $k\ge 2$, Equations \eqref{eqn:Z1}, \eqref{eqn:Z2} simplify to the pleasant formulas $Z_1=1$, $Z_2=-e^{\frac{-2\pi i}{3k}}$.
In particular, plugging into Equation \eqref{eqn:Uz} and simplifying gives that
$$U(Z_1)=\frac{3}{2}\left(\mathrm{Li}_2\left(e^{\frac{2\pi i}{3}+\frac{2\pi i}{3k}}\right)-\mathrm{Li}_2\left(-e^{\frac{\pi i}{3}+\frac{2\pi i}{3k}}\right)\right),\quad U(Z_2)=\frac{3}{2}\left(\mathrm{Li}_2\left(-e^{\frac{2\pi i}{3}}\right)-\mathrm{Li}_2\left(e^{\frac{\pi i}{3}}\right)\right).$$
We therefore compute from Equation \eqref{eqn:volumetetra} that the volume of this tetrahedron is $\frac{1}{2}\Im\bigl(U(Z_1) - U(Z_2)\bigr)$ and, as $M$ is comprised of $2k$ such tetrahedra, we find that
\begin{align}
\mathrm{vol}(M_k)&=\frac{3k}{2}\Im\left( \mathrm{Li}_2\left(e^{\frac{2\pi i}{3}+\frac{2\pi i}{3k}}\right)+\mathrm{Li}_2\left(e^{\frac{\pi i}{3}}\right)-\mathrm{Li}_2\left(-e^{\frac{2\pi i}{3}}\right)-\mathrm{Li}_2\left(-e^{\frac{\pi i}{3}+\frac{2\pi i}{3k}}\right)\right),\nonumber\\
&=\frac{3k}{2}\left( D\left(e^{\frac{2\pi i}{3}+\frac{2\pi i}{3k}}\right)+2D\left(e^{\frac{\pi i}{3}}\right)+D\left(e^{\frac{2\pi i}{3}-\frac{2\pi i}{3k}}\right)\right),\label{eqn:volumeMkk}
\end{align}
where $D(-)$ is the Bloch--Wigner function defined on roots of unity by
$$D(e^{i\theta})=\Im(\mathrm{Li}_2(e^{i\theta}))=\sum_{n=1}^\infty\frac{\sin(n\theta)}{n^2},$$
and the simplification in the second line comes from routine properties of sines.
This gives the requisite volume.

To see that $\mathrm{vol}(M_k)$ is linear in $k$, one notes (from \cite[\S 3]{Zagier} for instance) that 
$$1\le D\left(e^{\frac{2\pi i}{3}+\frac{2\pi i}{3k}}\right)+D\left(e^{\frac{2\pi i}{3}-\frac{2\pi i}{3k}}\right),\quad\quad 1\le D\left(e^{\frac{\pi i}{3}}\right),$$
and therefore 
$$\frac{9k}{2}\le \mathrm{vol}(M_k)\le \pi^2k,$$
as required.


\subsection{Computational examples for \texorpdfstring{$k$}{k} small}


The unique manifold $M_2$ in $\M_{2,2}$ has been identified computationally (see \cite[\S 7.6]{Anuthesis} or \cite{fmp-small}), as has its double $DM_2$.
One can use Snappy \cite{snappy} to compute its volume, homology, and isometry group using the following code:
\begin{verbatim}
DM=Manifold(`tLwvLLLwAMPQQkbbhhjilkmmmsporqqrsspuappuoowsfpabuuvvbb') 
DM.volume()
DM.homology()
DM.symmetry_group()
\end{verbatim}
Consequently, we find that $DM_2$ is a four-cusped hyperbolic manifold with an embedded totally geodesic surface of genus $2$, which has volume approximately $18.2689489153$, nine times the volume of the figure-eight knot complement, and homology given by $\Z^4\oplus\Z/3\Z$, with each copy of $\Z$ being generated by the meridian of a cusp.
Note that the volume of $DM_2$ is exactly twice the value of Equation \eqref{eqn:volumeMkk} when evaluated at $k=2$, as one expects from doubling.
Moreover, Snappy computes that $\Isom(DM_2)\cong D_6\times\Z/2\Z$, which perfectly aligns with the statement of Theorem \ref{thm-isometries-kk}.
Indeed, the $D_6$-factor of this isometry group is generated by performing the isometries from Theorem \ref{thm-isometries-kk} on each of the two copies of $M_2$ in the double $DM_2$ and the $\Z/2\Z$-factor is the isometry which interchanges these two factors.
Moreover, one can compute the arithmetic invariants of $DM_2$ using Snap in Snappy via the following code:
\begin{verbatim}
from snappynt.ManifoldNT import ManifoldNT as snappyNT
DM=snappyNT(`tLwvLLLwAMPQQkbbhhjilkmmmsporqqrsspuappuoowsfpabuuvvbb')
DM.compute_arithmetic_invariants()
DM.p_arith()
\end{verbatim}
The resulting computation shows that $DM_2$ has invariant trace field $\Q(\sqrt{-3})$, invariant quaternion algebra given by the matrix algebra $\mathrm{Mat}_2(\Q(\sqrt{-3}))$, and non-integral traces.
This aligns with the statement of Theorem \ref{thm:kevenquasiarithmetic} and confirms computationally that $DM_2$ is a non-arithmetic, quasi-arithmetic manifold.

Note that Corollary \ref{cor:kplusonecount} shows that $|\mathcal{M}_{2,1}|=1$, we denote the unique manifold in this set by $M_1$ and its double by $DM_1$.
Then Theorem \ref{thm-dehnfill} shows that performing any of the six Dehn fillings listed therein on a cusp of $M_2$ produces $M_1$. 
The operations of Dehn fillings on a cusp and doubling commute so one may construct $DM_1$ by performing the same Dehn filling on two of the four cusps of $DM_2$, which are interchanged by the aforementioned $\Z/2\Z$-isometry. 
This can be seen computationally in Snappy by the following commands.
\begin{verbatim}
DM.dehn_fill((2,3),2)
DM.dehn_fill((3,2),3)
\end{verbatim}
The resulting manifold $DM_1$ is a two-cusped hyperbolic manifold with an embedded totally geodesic surface of genus $2$, which has volume approximately $15.5952737606$ and homology given by $\Z^3$.
Once again, Snappy computes that $\Isom(DM_1)\cong \Z/2\Z$, which aligns with the second part of Corollary \ref{cor:kplusonecount}, since the non-trivial element in this group is the isometry switching the two copies of $M_1$ in $DM_1$.
Performing the same computational routine to compute the arithmetic invariants as above, one also learns that the invariant trace field is a degree $6$ field with $3$ complex places and minimal polynomial $x^6 - x^5 + 4x^4 - 4x^3 + 5x^2 - 3x + 1$.
Moreover, the invariant quaternion algebra is again a matrix algebra over this field, as must be the case. 


\bibliographystyle{amsalpha}
\bibliography{bib}

\end{document}